\documentclass[a4paper,11pt]{amsart}

\usepackage{amssymb,amsmath}
\usepackage{graphicx}
\usepackage[all]{xy}

\title[A model category for the homotopy theory of concurrency]{A model category for the homotopy theory of concurrency}
\author[P. Gaucher]{Philippe Gaucher}
\email{gaucher@pps.jussieu.fr}
\urladdr{http://www.pps.jussieu.fr/{\~{}}gaucher/}
\address{Preuves Programmes et Syst{\`e}mes\\ Universit{\'e} Paris 7--Denis Diderot\\
Case 7014\\2 Place Jussieu\\ 75251 PARIS Cedex 05\\ France}

\newcommand{\C}{\mathcal{C}}
\newcommand{\D}{\mathcal{D}}
\newcommand{\Z}{\mathbb{Z}}

\newcommand{\R}{\mathbb{R}}
\newcommand{\de}{\partial}
\newcommand{\p}\times
\renewcommand{\vec}{\overrightarrow}
\renewcommand{\P}{\mathbb{P}}

\newcommand{\be}{\begin{equation}}
\newcommand{\ee}{\end{equation}}
\newcommand{\bea}{\begin{eqnarray}}
\newcommand{\eea}{\end{eqnarray}}
\newcommand{\beas}{\begin{eqnarray*}}
\newcommand{\eeas}{\end{eqnarray*}}

\newtheorem{thm}{Theorem}[section]
\newtheorem{prop}[thm]{Proposition}

\newtheorem{cor}[thm]{Corollary}
\newtheorem{rem}[thm]{Remark}
\newtheorem{qu}[thm]{Question}
\newtheorem{defn}[thm]{Definition}

\newtheorem{nota}[thm]{Notation}
\newcommand{\bd}{\begin{defn}}
\newcommand{\ed}{\end{defn}}
\newcommand{\bcd}{\begin{defn}}
\newcommand{\ecd}{\end{defn}}
\newcommand{\bex}{\begin{exmp}}
\newcommand{\eex}{\end{exmp}}
\newcommand{\bp}{\begin{prop}}
\newcommand{\ep}{\end{prop}}
\newcommand{\bth}{\begin{thm}}
\renewcommand{\eth}{\end{thm}}
\newcommand{\br}{\begin{rem}}
\newcommand{\er}{\end{rem}}
\newcommand{\bpf}{\begin{proof}}
\newcommand{\epf}{\end{proof}}

\newcommand{\fl}[1]{\ar@{->}[l]_{#1}}
\newcommand{\fr}[1]{\ar@{->}[r]^-{#1}}
\newcommand{\fd}[1]{\ar@{->}[d]_{#1}}
\newcommand{\fu}[1]{\ar@{->}[u]^{#1}}
\newcommand{\f}[2]{\ar@{->}[#1]|{#2}}
\newcommand{\ff}[2]{\ar@2{->}[#1]|{#2}}
\newcommand{\frr}[1]{\ar@{->}[rr]^{#1}}

\renewcommand{\top}{{\mathbf{Top}}}

\newcommand{\iso}{\cong}

\newcommand{\ot}{\otimes}

\newcommand{\vI}{\vec{I}}

\renewcommand{\leq}{\leqslant}
\renewcommand{\geq}{\geqslant}

\newcommand{\cattopn}{{\brm{1Cat^{top}_1}}}
\newcommand{\tcattopn}{{\brm{1CAT^{top}_1}}}

\newcommand{\homcat}{[\cattopn]}

\def\cartesien{%
  \ar@{-}[]+R+<6pt,-2pt>;[]+RD+<6pt,-6pt>%
  \ar@{-}[]+D+<2pt,-6pt>;[]+RD+<6pt,-6pt>%
}
\def\cocartesien{%
  \ar@{-}[]+L+<-6pt,+2pt>;[]+LU+<-6pt,+6pt>%
  \ar@{-}[]+U+<-2pt,+6pt>;[]+LU+<-6pt,+6pt>%
}

\newcommand{\brm}[1]{\rm{\mathbf{#1}}}

\renewcommand{\top}{{\brm{Top}}}
\newcommand{\gltop}{{\brm{glTop}}}

\newcommand{\dtop}{{\brm{Flow}}}

\newcommand{\set}{{\brm{Set}}}
\newcommand{\tdtop}{{\brm{FLOW}}}
\newcommand{\ttop}{{\brm{TOP}}}

\newcommand{\sets}{{\brm{Set}}}

\newcommand{\glob}{{\rm{Glob}}}

\newcommand{\liminj}{\varinjlim}
\newcommand{\limproj}{\varprojlim}

\makeatletter

\def\varholim@#1#2{%
  \vtop{\m@th\ialign{##\cr
    \hfil$#1\operator@font holim$\hfil\cr
    \noalign{\nointerlineskip\kern1.5\ex@}#2\cr
    \noalign{\nointerlineskip\kern-\ex@}\cr}}%
}
\def\holimproj{%
  \mathop{\mathpalette\varholim@{\leftarrowfill@\textstyle}}\nmlimits@
}
\def\holiminj{%
  \mathop{\mathpalette\varholim@{\rightarrowfill@\textstyle}}\nmlimits@
}

\newskip\@bigflushglue \@bigflushglue = -100pt plus 1fil

\def\bigcentering{\let\\\@centercr\rightskip\@bigflushglue%
\leftskip\@bigflushglue
\parindent\z@\parfillskip\z@skip}

\makeatother

\DeclareMathOperator{\id}{Id} \DeclareMathOperator{\Id}{Id}
\DeclareMathOperator{\card}{card}

\newcommand{\diag}{\mathcal{D}}

\begin{document}

\begin{abstract}
We construct a cofibrantly generated model structure on the category
of flows such that any flow is fibrant and such that two cofibrant
flows are homotopy equivalent for this model structure if and only if
they are S-homotopy equivalent. This result provides an interpretation
of the notion of S-homotopy equivalence in the framework of model
categories.
\end{abstract}

\maketitle

\tableofcontents

\section{Geometric models of  concurrency}

Algebraic topological models have been used now for some years in
concurrency theory (concurrent database systems and fault-tolerant
distributed systems as well) \cite{survol}. The earlier models, {\em
progress graph} (see \cite{CoElSh71} for instance) have actually
appeared in operating systems theory, in particular for describing the
problem of ``deadly embrace'' (as E. W.  Dijkstra originally
put it in \cite{EWDCooperating}, now more usually called deadlock) in
``multiprogramming systems''. They are used by J. Gunawardena in
\cite{Gunawardena1} as an example of the use
of homotopy theory in concurrency theory. Later V. Pratt introduced
another geometric approach using \textit{strict globular
$\omega$-categories} in \cite{Pratt}. Some of his ideas would be
developed in an homological manner in E. Goubault's PhD \cite{HDA},
using bicomplexes of modules.  The $\omega$-categorical point of view
would be developed by the author mainly in \cite{Gau} \cite{Coin}
\cite{sglob} \cite{fibrantcoin} using the equivalence of
categories between the category of strict globular $\omega$-categories
and that of strict cubical $\omega$-categories
\cite{math.CT/0007009}. The mathematical works of R. Brown \textit{et
al.} \cite{Brown_cube} \cite{RBPJHColimit} and of R.  Street
\cite{oriental} play an important role in this approach.

The $\omega$-categorical approach also allowed to understand how to
deform \textit{higher dimensional automata} (HDA) modeled by
$\omega$-categories without changing their compu\-ter-scientific
properties (deadlocks, unreachable states, schedules of execution,
final and initial points, serializability). The notions of
\textit{spatial deformation} and of \textit{temporal deformation} of
HDA are indeed introduced in \cite{ConcuToAlgTopo} in an informal way.

Another algebraic topological approach of concurrency is that of
\textit{local po-space} introduced by L. Fajstrup, E. Goubault and
M. Raussen.  A local po-space is a gluing of topological spaces which
are equipped with a closed partial ordering representing the time
flow. They are used as a formalization of higher dimensional automata
which model concurrent systems in computer science. Some algorithms of
deadlock detection in PV diagrams have been studied within this
framework \cite{LFEGMRDetecting}.

The notion behind all these geometric approaches is the one of
\textit{precubical set}. Roughly speaking, a $n$-dimensional cube
$[0,1]^n$ represents the concurrent execution of $n$ independant
processes. A precubical set is a family of sets $(K_n)_{n\geq 0}$ (the
elements of $K_n$ being called the $n$-dimensional cubes) together
with face operators $\partial_i^\alpha:K_{n+1}\longrightarrow K_{n}$
for $1\leq i\leq n$ and with $\alpha\in\{-,+\}$ satisfying
$\partial_i^\alpha\partial_j^\beta=\partial_{j-1}^\beta\partial_i^\alpha$
for $i<j$. These face operators encode how the $n$-cubes are located
with respect to one another in the precubical set. The prefix ``pre''
means that there are no degeneracy maps at all in the data.
R. Cridlig presents in \cite{cridlig96implementing} an implementation
with CaML of the semantics of a real concurrent language in terms of
\textit{precubical sets}, demonstrating the relevance of this
approach. Since this category is sufficient to model HDA, why not
deal directly with precubical sets ?  Because the category of
precubical sets is too poorly structured.  For instance there are not
enough morphisms to model temporal deformations (see also the
introduction of \cite{Coin} for some further closely related reasons).

In \cite{diCW}, some particular cases of local po-spaces are
introduced by E. Goubault and the author: \textit{the globular
CW-complexes}. The corresponding category is big enough to model all
HDA. Moreover the notion of spatial and temporal deformations can be
modeled within this category. It became possible to give a precise
mathematical definition of two globular CW-complexes to be
\textit{S-homotopy equivalent} and \textit{T-homotopy equivalent} (S
for space and T for time !). By localizing with respect to the
S-homotopy and T-homotopy equivalences, one obtains a new category,
that of \textit{dihomotopy types}, whose isomorphism classes are
globular CW-complexes having the same computer scientific
properties. It then became possible to study concurrency using only
this quotient category of dihomotopy types.

Not only globular complexes allow to model dihomotopy, but they also
allow to take out pathological situations appearing in the local
po-space framework and which are meaningless from a computer
scientific viewpoint. For example, the rational numbers $\mathbb{Q}$
equipped with the usual ordering is a local po-space and the total
disconnectedness of $\mathbb{Q}$ means nothing in this geometric
approach of concurrency.

The purpose of this paper is the introduction of a new category, the
category of \textit{flows}, in which it will be possible to embed the
category of globular CW-complexes and in which it will be possible to
define both the class of S-homotopy and T-homotopy equivalences. Due
to the length of this work, the construction and the study of the
functor from the category of globular CW-complexes to that of flows is
postponed to another paper.

\begin{figure}
\begin{center}
\includegraphics[width=12cm]{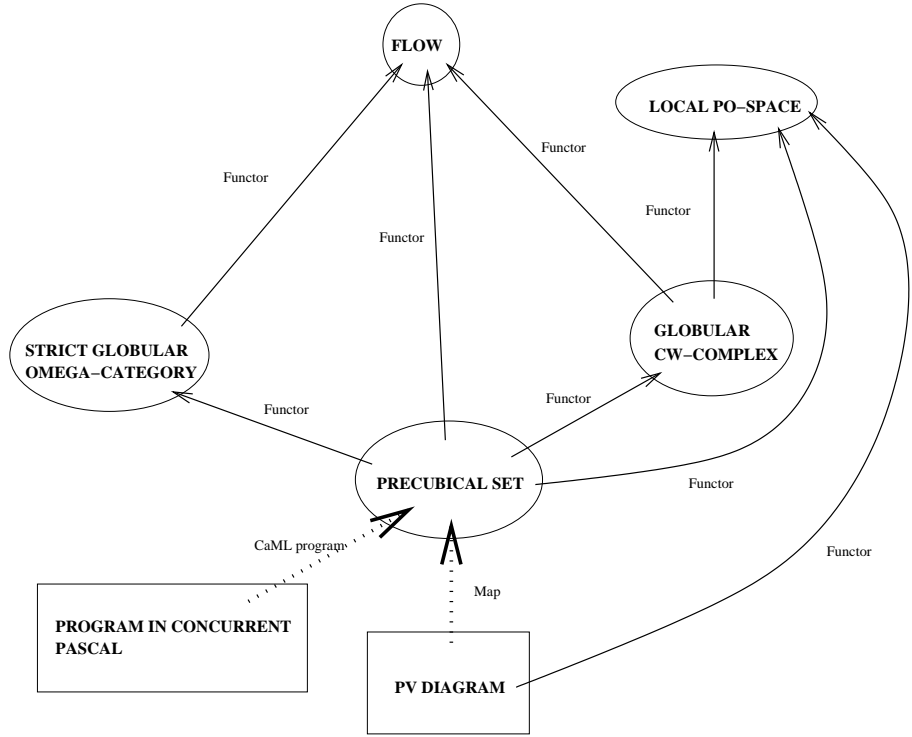}
\end{center}
\caption{Comparison of geometric models of HDA} \label{form}
\end{figure}

Figure~\ref{form} is a recapitulation of the geometric models of
concurrency, including the one presented in this paper.

\section{Outline of the paper}

Section~\ref{def} defines the category of flows $\dtop$ after a short
introduction about compactly generated topological spaces.  It is
proved that $\dtop$ is complete and cocomplete. Several particular and
important examples of flows are also introduced.  Section~\ref{com} is
devoted to proving that for any flow $Y$, the functor $\tdtop(-,Y)$
from the opposite of the category of flows to that of topological
spaces commutes with all limits where $\tdtop(X,Y)$ is the set of
morphisms of flows from $X$ to $Y$ endowed with the Kelleyfication of
the compact-open topology.  This fact will be of crucial importance in
several places of the paper. This result turns out to be difficult to
establish since the underlying topological space of a colimit of flows
is in general not isomorphic to the colimit of the underlying
topological spaces. This result actually requires the introduction of
the category of \textit{non-contracting topological $1$-categories}
and of a closed monoidal structure on it. Section~\ref{canonique}
shows that any flow is a canonical colimits of globes and points. This
is a technical lemma which is also of importance for several proofs of
this paper. Section~\ref{sectionS} defines the class of S-homotopy
equivalences in the category of flows. The associated cylinder functor
is constructed. Section~\ref{boxex} is devoted to an explicit
description of $U\boxtimes X$ for a given topological space $U$ and a
given flow $X$. Section~\ref{SDEP} describes a class of morphisms of
flows (the ones satisfying the S-homotopy extension property) which
are closed by pushouts and which contains useful examples as the
inclusion $\glob(\de Z)\longrightarrow
\glob(Z)$ where $(Z,\de Z)$ is a NDR pair of topological spaces.
The main result of Section~\ref{secincl} is that any morphism of flows
satisfying the S-homotopy extension property induces a closed
inclusion of topological spaces between the path spaces. This allows
us to prove in Section~\ref{small} that the domains of the generating
cofibrations and of the generating trivial cofibrations of the model
structure are small relatively to the future class of cofibrations of
the model structure. Section~\ref{small} is therefore the beginning of
the construction of the model structure. Section~\ref{remindermodel}
recalls some well-known facts about cofibrantly generated model
categories.  Section~\ref{chafib} characterizes the fibrations of this
model structure. Section~\ref{necessaire} explains why it is necessary
to add to the set of generating cofibrations the morphisms of flows
$C:\varnothing\longrightarrow \{0\}$ and $R:\{0,1\}\longrightarrow
\{0\}$. Section~\ref{expush} provides an explicit calculation of the pushout
of a morphism of flows of the form $\glob(\de Z)\longrightarrow
\glob(Z)$. This will be used in Section~\ref{fin1}. The main result of
Section~\ref{expush} is that if $\de Z\longrightarrow Z$ is an
inclusion of a deformation retract, then any morphism of flows which
is a pushout of $\glob(\de Z)\longrightarrow \glob(Z)$ induces a weak
homotopy equivalence between path spaces.  Section~\ref{fin1} and
Section~\ref{fin2} conclude the construction of the model structure
recapitulated in Section~\ref{modelflow}. Section~\ref{comparison}
checks that two cofibrant-fibrant flows are homotopy equivalent for
this model structure if and only if they are S-homotopy equivalent.

\section{Warning}

This paper is the first part of a work which aims at introducing a
convenient categorical setting for the homotopy theory of concurrency.
This part is focused on the category of flows itself, its basic
properties, the notion of S-homotopy equivalence, weak or not, and the
model structure. The relation between the category of globular
CW-complexes and the one of flows is explored in \cite{model2}. A
detailed abstract (in French) of this work can be found in
\cite{pgnote1} and \cite{pgnote2}.

\section{The category of flows}\label{def}

\subsection{Preliminaries about the compactly generated topological spaces}

This section is a survey about compactly generated spaces which gives
enough references for the reader not familiar with this subject. Cf.
\cite{MR90k:54001}, \cite{MR2000h:55002} and the appendix of \cite{Ref_wH}.

By a \textit{compact space}, we mean a compact Hausdorff topological space.
Let $\mathcal{T}$ be the category of general topological spaces with the
continuous maps as morphisms.

\bd A continuous map $f:A\longrightarrow B$ is an {\rm inclusion
of spaces} if $f$ is one-to-one and if the canonical set map
\[\top(Z,A)\longrightarrow \{g\in \top(Z,B),g(Z)\subset f(A)\}\]
induced by the mapping $g\mapsto f\circ g$ is a bijection of sets. In
other terms, a continuous map $f:A\longrightarrow B$ is an inclusion
of spaces if for any set map $g:Z\longrightarrow A$ such that $f\circ
g$ is continuous, then $g$ is continuous.
\ed

\bd A continuous map $f:A\longrightarrow B$ is {\rm closed} if for any
closed subset $F$ of $A$, the subset $f(F)$ is closed in $B$. \ed

\bd A {\rm quotient map} is a continuous map $f:X\longrightarrow Y$
which is onto and such that $U\subset Y$ is open if and only if
$f^{-1}(U)$ is open in $X$. In other term, $Y$ is given with the final
topology associated to $f$. \ed

\bd A {\rm $k$-space} $X$ is a topological space such that for any
continuous map $f:K\longrightarrow X$ with $K$ compact, $U\subset X$
is open (resp.  closed) if and only if $f^{-1}(U)$ is open
(resp. closed) in $K$.  The corresponding category with the continuous
maps as morphisms is denoted by $k\top$. \ed

A topological space $X$ is a $k$-space if and only if there
exists a disjoint sum of compacts $\bigoplus_{i\in I} K_i$ and a
quotient map $\bigoplus_{i\in I} K_i\longrightarrow X$
\cite{MR90k:54001}. The inclusion functor $k\top\longrightarrow
\mathcal{T}$ has a right adjoint and a left inverse
$k:\mathcal{T}\longrightarrow k\top$ which is called the
\textit{Kelleyfication} functor.  The category $k\top$ is
complete and cocomplete where colimits are taken in $\mathcal{T}$ and
limits are taken by applying $k$ to the limit in $\mathcal{T}$
\cite{MR35:970} \cite{MR1712872}.  The identity map
$k\left(X\right)\longrightarrow X$ is continuous because the topology
of $k\left(X\right)$ contains more opens than the topology of $X$.

\bd A topological space $X$ is {\rm weak Hausdorff} if and only if
for any continuous map $f:K\longrightarrow X$ with $K$ compact, the
subspace $f(K)$ is closed in $X$. \ed

If $X$ is a $k$-space, then $X$ is weak Hausdorff if and only if
its diagonal $\Delta X=\{(x,x)\in X\p X\}$ is a closed subspace
of $X\p X$, the latter product being taken in $k\top$
\cite{MR40:4946}. If $X$ is a weak Hausdorff topological space,
then $k(X)$ is still weak Hausdorff.

If $X$ is a weak Hausdorff topological space, then $X$ is a $k$-space
if and only if $X\iso \liminj_{K\subset X} K$ as topological space
where $K$ runs over the set of compact subspaces of $X$: a subset $F$
of $k\left(X\right)$ is closed (resp. open) if and only if for any
compact $C$ of $X$, $F\cap C$ is a closed (resp. open) subspace of
$X$.

\bd A {\rm compactly generated topological space} is by definition
a weak Hausdorff $k$-space. The corresponding category with the
continuous maps as morphims is denoted by $\top$. \ed

Let $wH$ be the category of weak Hausdorff topological spaces.
Generally colimits in $wH$ do not coincide with colimits in
$\mathcal{T}$. But

\bp\cite{Ref_wH} A transfinite composition of injections and
push\-outs of closed inclusions of compactly generated topological
spaces is still weak Hausdorff (and therefore a compactly generated
topological space). \ep

\bp\cite{MR40:4946} \cite{MR1712872} The inclusion functor
$wH\longrightarrow \mathcal{T}$ has a left adjoint $H$. If $X$ is a
$k$-space and if $\mathcal{R}$ is an equivalence relation, then
$H(X/\mathcal{R})$ is equal to $X/\overline{\mathcal{R}}$ where the
topological closure $\overline{\mathcal{R}}$ of $\mathcal{R}$ is
defined as the intersection of all equivalence relations containing
$\mathcal{R}$ and whose graph is closed in $X\p X$. In particular, if
the graph of $\mathcal{R}$ is closed in $X\p X$, then $X/\mathcal{R}$
is weak Hausdorff. \ep

\bp\label{limitset} \cite{MR35:970} \cite{MR1712872} If $i\mapsto X_i$
is any small diagram in $\top$, then the limit in $\top$ coincides
with the Kelleyfication of the limit in $\mathcal{T}$ and with the
Kelleyfication of the limit in $wH$. Moreover the underlying set of
this limit coincides with the limit in the category of sets of the
underlying sets of the $X_i$. \ep

If $X$ is a weak Hausdorff topological space, then a subset $Y$
of $X$ equipped with the relative topology is weak Hausdorff as
well. If $X$ is a compactly generated topological space, then a
subset $Y$ of $X$ equipped with the relative topology is then
weak Hausdorff. But it is not necessarily a $k$-space. To get
back a $k$-space, it is necessary to consider the Kelleyfication
$k(Y_r)$ of $Y_r$ ($Y$ equipped with the relative topology).

\bp\label{lim}\cite{MR35:970} \cite{MR1712872} Let us denote by
$\ttop\left(X,-\right)$ the right adjoint of the functor $-\p
X:\top\longrightarrow \top$. Then
\begin{enumerate}
\item If $Cop\left(X,Y\right)$ is the set $\top\left(X,Y\right)$ equipped with the
compact-open\index{compact-open topology} topology (i.e. a basis
of opens is given by the sets
\[N\left(C,U\right):=\{f\in\top\left(X,Y\right),f\left(C\right)\subset
U\}\] where $C$ is any compact subset of $X$ and U any open subset
of $Y$), then there is a natural bijection
$\ttop\left(X,Y\right)\iso k\left(Cop\left(X,Y\right)\right)$.
\item There is a natural isomorphism of topological spaces
\[\ttop\left(X\p Y,Z\right)\iso
\ttop\left(X,\ttop\left(Y,Z\right)\right).\]
\item There are natural isomorphisms of topological spaces
\[\ttop\left(\liminj_i X_i,Y\right)\iso \limproj_i \ttop\left(X_i,Y\right)\] and
\[\ttop\left(X,\limproj_i Y_i\right)\iso \limproj_i \ttop\left(X,Y_i\right).\]
\end{enumerate}
\ep

Similar results can be found in \cite{MR49:11475} \cite{MR45:9323}
with slightly bigger categories of topological spaces than the one we
are using in this paper.

In the sequel, all topological spaces will be supposed to be compactly
generated (so in particular weak Hausdorff). In particular all binary
products will be considered within this category.

\subsection{Definition of a flow}

\bd A {\rm flow} $X$ consists of a topological space $\P X$, a
discrete space $X^0$, two continuous maps $s$ and $t$ from $\P X$ to
$X^0$ and a continuous and associative map \[*:\{(x,y)\in \P X\p \P X;
t(x)=s(y)\}\longrightarrow \P X\] such that $s(x*y)=s(x)$ and
$t(x*y)=t(y)$.  A morphism of flows $f:X\longrightarrow Y$ consists of
a set map $f^0:X^0\longrightarrow Y^0$ together with a continuous map
$\P f:\P X\longrightarrow \P Y$ such that $f(s(x))=s(f(x))$,
$f(t(x))=t(f(x))$ and $f(x*y)=f(x)*f(y)$. The corresponding category
will be denoted by $\dtop$. \ed

The continuous map $s:\P X\longrightarrow X^0$ is called the
\textit{source map}. The continuous map $t:\P X\longrightarrow X^0$
is called the \textit{target map}. One can canonically extend these
two maps to the whole underlying topological space $X^0\sqcup \P X$ of
$X$ by setting $s\left(x\right)=x$ and $t\left(x\right)=x$ for $x\in
X^0$.

The discrete topological space $X^0$ is called the
\textit{$0$-skeleton\index{$0$-skeleton}} of $X$. The
$0$-dimensional elements of $X$ are also called \textit{states} or
\textit{constant execution paths}.

The elements of $\P X$ are called \textit{non constant execution
paths}.  If $\gamma_1$ and $\gamma_2$ are two non-constant
execution paths, then $\gamma_1 *\gamma_2$ is called the
\textit{concatenation} or the \textit{composition} of $\gamma_1$
and $\gamma_2$. For $\gamma\in \P X$, $s\left(\gamma\right)$ is called
the \textit{beginning} of $\gamma$ and $t\left(\gamma\right)$ the
\textit{ending} of $\gamma$.

\begin{nota} For $\alpha,\beta\in X^0$, let $\P_{\alpha,\beta}X$ be the
subspace of $\P X$ equipped the Kelleyfication of the relative
topology consisting of the non-execution path\index{execution path}s
of $X$ with beginning $\alpha$ and with ending $\beta$. \end{nota}

\bd Let $X$ be a flow. A point $\alpha$ of $X^0$ such that there is not 
any non-constant execution path $\gamma$ with
$t\left(\gamma\right)=\alpha$ (resp.  $s\left(\gamma\right)=\alpha$)
is called {\rm an initial state} (resp. {\rm a final state}). \ed

\subsection{The globe of a topological space}

As in \cite{diCW}, but here for the framework of flows, we are going
to introduce the notion of \textit{globe} of a topological space. It
will be important both for computer scientific and purely mathematical
reasons.

\begin{figure}
\begin{center}
\includegraphics[width=7cm]{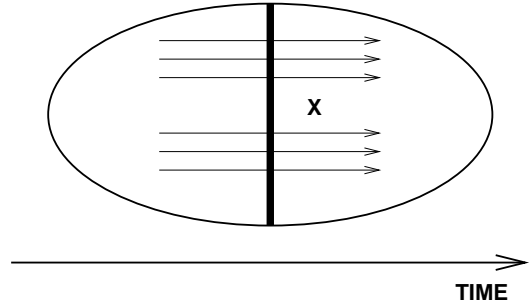}
\end{center}
\caption{Symbolic representation of
$\glob(X)$ for some topological space $X$} \label{exglob}
\end{figure}

For $X$ a topological space, let $\glob\left(X\right)$ be the flow
defined by
\[\glob\left(X\right)^0=\{0,1\}\hbox{ and }\P \glob\left(X\right)=X\]
with $s=0$ and $t=1$ (cf. Figure~\ref{exglob}). The Glob mapping
induces a canonical functor from the category $\top$ of topological
spaces to the category $\dtop$ of flows.

As a particular case of globe is that of a singleton. One obtains the
\textit{directed segment} $\vI$. It is defined as follows:
$\vI^0=\{0,1\}$, $\P\vI=\{[0,1]\}$, $s\left([0,1]\right)=0$ and
$t\left([0,1]\right)=1$.

If $Z_1,\dots,Z_p$ are $p$ topological spaces with $p\geq 2$, the
flow
\[\glob(Z_1)*\glob(Z_2)*\dots *\glob(Z_p)\]
is the flow obtained by identifying the final state of $\glob(Z_i)$
with the initial state of $\glob(Z_{i+1})$ for $1\leq i\leq p-1$.

\begin{nota} If $X$ and $Y$ are two flows, let us denote by $\tdtop(X,Y)$
the space of morphisms of flows $\dtop(X,Y)$ equipped with the
Kelleyfication of the compact-open topology.
\end{nota}

\bp \label{can} Let $X$ be a flow. Then there is a natural
homeomorphism $\P X\iso \tdtop\left(\vI,X\right)$. \ep

\bpf If we have an element $u$ of $\P X$, consider the morphism
of flows $F_\gamma$ defined by $F_\gamma\left(0\right)=s (u)$,
$F_\gamma\left(1\right)=t (u)$ and $F_\gamma\left([0,1]\right)=u$.
And reciprocally a morphism $F\in \dtop\left(\vI,X\right)$ can be
mapped on an element of $\P X$ by $F\mapsto F\left([0,1]\right)$.
Hence the bijection between the underlying sets. This bijection is an
homemorphism since for any topological space $Z$, one has the
homeomorphism $\ttop(\{0\},Z)\iso Z$.  \epf

\subsection{Higher dimensional automaton and flow}

This example is borrowed from \cite{diCW}. An example of {\em progress
graph}, that is of higher dimensional automaton, is modeled here as a
flow.

The basic idea is to give a description of what can happen when
several processes are modifying shared resources. Given a shared
resource $a$, we see it as its associated semaphore that rules its
behaviour with respect to processes. For instance, if $a$ is an
ordinary shared variable, it is customary to use its semaphore to
ensure that only one process at a time can write on it (this is mutual
exclusion). A semaphore is nothing but a register which counts the
number of times a shared object can still be accessed by processes. In
the case of usual shared variables, this register is initialized with
value 1, processes trying to access (read or write) on the
corresponding variable compete in order to get it first, then the
semaphore value is decreased: we say that the semaphore has been
locked\footnote{Of course this operation must be done ``atomically'',
meaning that the semaphore itself must be handled in a mutually
exclusive manner: this is done at the hardware level.} by the process.
When it is equal to zero, all processes trying to access this
semaphore are blocked, waiting for the process which holds the lock to
relinquish it, typically when it has finished reading or writing on
the corresponding variable: the value of the semaphore is then
increased.

When the semaphores are initialized with value one, meaning that they
are associated with shared variables accessed in a mutually exclusive
manner, they are called binary semaphores. When a shared data
(identified with its semaphore) can be accessed by one or more
processes, meaning that the corresponding semaphore has been
initialized with a value greater than one, it is called a counting
semaphore.

Given $n$ deterministic sequential processes $Q_1,\ldots,Q_n$,
abstracted as a sequence of locks and unlocks on (semaphores
associated with) shared objects, \[Q_i=R^1 a_i^1.R^2 a_i^2 \cdots
R^{n_i} a_i^{n_i}\] ($R^k$ being $P$ or $V$\footnote{Using E. W.
Dijkstra's notation $P$ and $V$ \cite{EWDCooperating} for respectively
acquiring and releasing a lock on a semaphore.}), there is a natural
way to understand the possible behaviours of their concurrent
execution, by associating to each process a coordinate line in
$\R^n$. The state of the system corresponds to a point in $\R^n$,
whose $i$th coordinate describes the state (or ``local time'') of the
$i$th processor.

\begin{figure}
\begin{center}
\includegraphics[width=7cm]{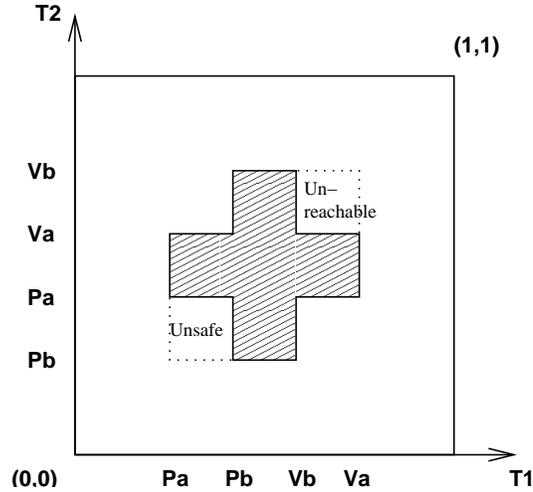}
\end{center}
\caption{Example of a progress graph}
\label{progress1}
\end{figure}

Consider a system with finitely many processes running altogether.  We
assume that each process starts at (local time) 0 and finishes at
(local time) 1; the $P$ and $V$ actions correspond to sequences of
real numbers between 0 and 1, which reflect the order of the $P$'s and
$V$'s. The initial state is $(0,\dots ,0)$ and the final state is
$(1,\dots ,1)$.  An example consisting of the two processes $T_1=P a.P
b.V b.Va$ and $T_2=P b.P a.V a.V b$ gives rise to the two dimensional
{\em progress graph\/} of Figure \ref{progress1}.

The shaded area represents states which are not allowed in any
execution path, since they correspond to mutual exclusion. Such states
constitute the {\em forbidden area}. An {\em execution path} is a path
from the initial state $(0,\ldots,0)$ to the final state
$(1,\ldots,1)$ avoiding the forbidden area and increasing in each
coordinate - time cannot run backwards. This entails that paths
reaching the states in the dashed square underneath the forbidden
region, marked ``unsafe'' are deemed to deadlock, i.e. they cannot
possibly reach the allowed terminal state which is $(1,1)$
here. Similarly, by reversing the direction of time, the states in the
square above the forbidden region, marked ``unreachable'', cannot be
reached from the initial state, which is $(0,0)$ here.  Also notice
that all terminating paths above the forbidden region are
``equivalent'' in some sense, given that they are all characterized by
the fact that $T_2$ gets $a$ and $b$ before $T_1$ (as far as resources
are concerned, we call this a {\em schedule}). Similarly, all paths
below the forbidden region are characterized by the fact that $T_1$
gets $a$ and $b$ before $T_2$ does.

We end up the paragraph with the Swiss Flag example of
Figure~\ref{progress1} described as a flow.

Let $n\geq 1$. Let $\mathbf{D}^n$ be the closed $n$-dimensional disk
defined by the set of points $\left(x_1,\dots,x_n\right)$ of $\R^n$
such that $x_1^2+\dots +x_n^2\leq 1$ endowed with the topology induced
by that of $\R^n$. Let $\mathbf{S}^{n-1}=\de \mathbf{D}^n$ be the
boundary of $\mathbf{D}^n$ for $n\geq 1$, that is the set of
$\left(x_1,\dots,x_n\right)\in \mathbf{D}^n$ such that $x_1^2+\dots
+x_n^2=1$. Notice that $\mathbf{S}^0$ is the discrete two-point
topological space $\{-1,+1\}$.  Let $\mathbf{D}^0$ be the one-point
topological space.  Let $\mathbf{S}^{-1}=\varnothing$ be the empty
set.

\begin{figure}
\begin{center}
\includegraphics[width=7cm]{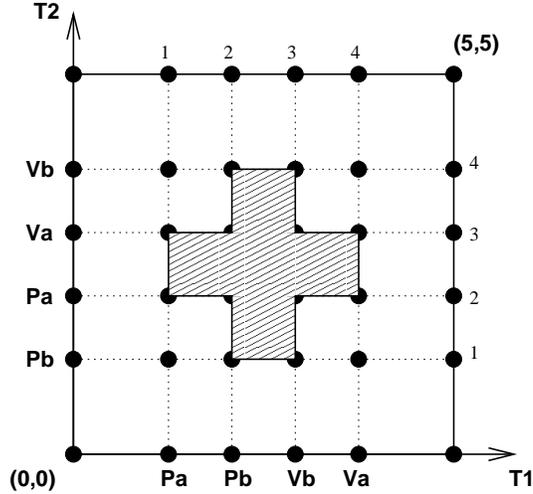}
\end{center}
\caption{Example of a flow}
\label{progress10}
\end{figure}

Consider the discrete set $SW^0=\{0,1,2,3,4,5\}\p \{0,1,2,3,4,5\}$.
Let
\beas
\mathcal{S}&=&\left\{((i,j),(i+1,j))\hbox{ for } (i,j)\in\{0,\dots,4\}\p \{0,\dots,5\}\right\}\\
&\cup& \left\{((i,j),(i,j+1))\hbox{ for } (i,j)\in\{0,\dots,5\}\p \{0,\dots,4\}\right\}\\
&\backslash & \left(
\{((2,2),(2,3)),((2,2),(3,2)), ((2,3),(3,3)),((3,2),(3,3))\}
\right)
\eeas
The flow $SW^1$ is obtained from $SW^0$ by attaching a copy of
$\glob(\mathbf{D}^0)$ to each pair $(x,y)\in \mathcal{S}$ with $x\in
SW^0$ identified with $0$ and $y\in SW^0$ identified with $1$.  The
flow $SW^2$ is obtained from $SW^1$ by attaching to each square
$((i,j),(i+1,j+1))$ except $(i,j)\in\{(2,1),(1,2),(2,2),(3,2),(2,3)\}$
a globular cell $\glob(\mathbf{D}^1)$ such that each execution path
$((i,j),(i+1,j),(i+1,j+1))$ and $((i,j),(i,j+1),(i+1,j+1))$ is
identified with one of the execution path of $\glob(\mathbf{S}^0)$
(there is not a unique choice to do that). Let $SW=SW^2$
(cf. Figure~\ref{progress10} where the bold dots represent the points
of the $0$-skeleton). The flow $SW$ represents the PV diagram of
Figure~\ref{progress10}.

\subsection{Limit and colimit in $\dtop$}

\bth\label{ssc}
\cite{MR96g:18001a} \cite{MR1712872} (Freyd's Adjoint Functor Theorem)
Let $A$ and $X$ be locally small categories. Assume that $A$ is
complete. Then a functor $G:A\longrightarrow X$ has a left adjoint if
and only if it preserves all limits and satisfies the following
``Solution Set Condition''. For each object $x\in X$, there is a set
of arrows $f_i:x\longrightarrow G a_i$ such that for every arrow
$h:x\longrightarrow G a$ can be written as a composite $h=Gt\circ f_i$
for some $i$ and some $t:a_i\longrightarrow a$. \eth

\bth \label{lim-colim} The category $\dtop$ is complete and
cocomplete. In particular, a terminal object is the flow $\mathbf{1}$
having the discrete set $\{0,u\}$ as underlying topological space with
$0$-skeleton\index{$0$-skeleton} $\{0\}$ and path space $\{u\}$. And an
initial object is the unique flow $\varnothing$ having the empty set
as underlying topological space.  \eth

\bpf Let $X:I\longrightarrow \dtop$ be a functor from a small
category $I$ to $\dtop$. Let $Y$ be the flow defined as follows:
\begin{enumerate}
\item The $0$-skeleton\index{$0$-skeleton} $Y^0$ of $Y$ is defined as being the
limit as sets $\limproj_{I} \left(X\left(i\right)^0\right)$ equipped
with the discrete topology.
\item Let $\alpha,\beta\in \limproj_{I} \left(X\left(i\right)^0\right)$ and let $\alpha_i$
(resp. $\beta_i$) be the image of $\alpha$ ($\beta$) in
$X\left(i\right)^0$. Then let $\P_{\alpha,\beta}Y:= \limproj_i
\P_{\alpha_i,\beta_i}X\left(i\right)$ where the limit is taken in
$\top$.
\item For $\alpha,\beta,\gamma\in \limproj_{I} \left(X\left(i\right)^0\right)$, let
$\alpha_i$ (resp. $\beta_i$, $\gamma_i$) be the image of $\alpha$
(resp. $\beta$, $\gamma$) in $X\left(i\right)^0$. Then the composition
map $*:\P_{\alpha,\beta}Y\p
\P_{\beta,\gamma}Y\longrightarrow \P_{\alpha,\gamma}Y$ is taken as
the limits of the $*_i:\P_{\alpha_i,\beta_i}X\left(i\right)\p
\P_{\beta_i,\gamma_i}X\left(i\right)\longrightarrow
\P_{\alpha_i,\gamma_i}X\left(i\right)$.
\end{enumerate}
One does obtain a flow which is the limit $\limproj_{i\in I}
X\left(i\right)$. To prove that $\dtop$ is cocomplete, it suffices to
prove that the constant diagram functor $\Delta_I$ from $\dtop$ to the
category $\dtop^I$ of diagrams in $\dtop$ over the small category $I$
has a left adjoint using Theorem~\ref{ssc}\index{solution set
condition}. The functor $\Delta_I$ commutes with
limit\index{limit}s. It suffices now to find a set of solutions.
Consider a diagram $D$ of $\dtop^I$. There is a class of solutions by
taking all morphisms $f:D\rightarrow \Delta_I Y$ for $Y$ running over
the category $\dtop$ and for $f$ running over the set of morphisms
from $D$ to $\Delta_I Y$. Then one can suppose that $Y$ is the subflow
generated by the image of $D$, so that the cardinal $\card(Y)$ of $Y$
satisfies $\card(Y)\leq
\aleph_0\p \card(D)$. Then it suffices to consider the set $\{Z_i,i\in I\}$
of isomorphism classes of flows whose underlying set is of cardinal
less than $\aleph_0\p \card(D)$. Then $\card(I)\leq 2^{\left(\aleph_0\p
\card(D)\right)^5}$. So $I$ is a set. Therefore
$\bigcup_{i\in I}\dtop^I\left(D,\Delta_I\left(Z_i\right)\right)$ is a
set as well. One has obtained a set of solutions. \epf

\section{Morphisms of flows and colimits}\label{com}

The aim of this section is the proof of the following theorem:

\bth (Theorem~\ref{commute}) Let $\tdtop\left(X,Y\right)$ be the
set of morphisms of flows from $X$ to $Y$ equipped with the
Kelleyfication of the compact-open\index{compact-open topology}
topology. Then the mapping \[\left(X,Y\right)\mapsto
\tdtop\left(X,Y\right)\] induces  a functor from $\dtop\p \dtop$
to $\top$ which is contravariant with respect to $X$ and
covariant with respect to $Y$. Moreover:
\begin{enumerate}
\item One has  the natural homeomorphism \[\tdtop\left(\liminj_i X_i,Y\right)\iso
\limproj_i \tdtop\left(X_i,Y\right)\] for any colimit $\liminj_i
X_i$ in $\dtop$.
\item One has the natural homeomorphism \[\tdtop\left(X,\limproj_i Y_i\right)\iso
\limproj_i\tdtop\left(X, Y_i\right)\] for any finite limit
$\limproj_i X_i$ in $\dtop$.
\end{enumerate}
\eth

\subsection{Non-contracting topological $1$-category}

\bd  A {\rm non-contracting topological $1$-category} $X$ is a
pair of compactly generated topological spaces $(X^0,\P X)$ together
with continuous maps $s$, $t$ and $*$ satisfying the same properties
as in the definition of flow except that $X^0$ is not necessarily
discrete.  The corresponding category is denoted by $\cattopn$. \ed

\bd A non-contracting topological $1$-category $X$ is {\rm
achronal} if $\P X=\varnothing$. \ed

\bth\label{rr} The category $\cattopn$ is complete and cocomplete.
The inclusion functor $\widetilde{\omega}:\dtop\longrightarrow
\cattopn$ preserves finite limits. \eth

\bpf Let $X:I\longrightarrow \cattopn$ be a functor from a small
category $I$ to $\cattopn$. Then consider the topological
$1$-category $Y$ defined as follows:
\begin{enumerate}
\item Let $Y^0:=\limproj_i X\left(i\right)^0$, the limit being taken in $\top$.
\item Let $\P Y:=\limproj_i \P X\left(i\right)$, the limit being taken in $\top$.
\item Let $Y=Y^0\sqcup \P Y$ equipped with the source map, target
map and composition law limits of the source maps, target maps
and composition laws of the $X\left(i\right)$.
\end{enumerate}
The $1$-category $Y$ is clearly the limit $\limproj X$ in
$\cattopn$. The cocompleteness of $\cattopn$ is then proved using the
``solution set condition'' recalled in
Theorem~\ref{ssc}\index{solution set condition} as in the proof of
Theorem~\ref{lim-colim}. A finite limit of discrete topological spaces
is discrete. So to be able to conclude that the functor
$\widetilde{\omega}$ preserves finite limits, it then suffices to
compare the construction of limits in $\dtop$ in the proof of
Theorem~\ref{lim-colim} and the construction of limits in $\cattopn$
in this proof. \epf

Using the above constructions, one sees that the
$0$-skeleton\index{$0$-skeleton} functor
\[(-)^0:\cattopn\longrightarrow \top\] does commute with any limit.
However the $0$-skeleton\index{$0$-skeleton} functor
$(-)^0:\dtop\longrightarrow \top$ only commutes with finite
limits. On the contrary, both $0$-skeleton\index{$0$-skeleton}
functors $(-)^0:\cattopn\longrightarrow \top$ and
$(-)^0:\dtop\longrightarrow \top$ do commute with any colimit.

The functor $\widetilde{\omega}:\dtop\longrightarrow \cattopn$ does
not preserve general limits. As counterexample, take the achronal
$1$-categories $\Z/p^n\Z$ equipped with the discrete topology and
consider the tower of maps $\Z/p^{n+1}\Z\longrightarrow \Z/p^{n}\Z$
defined by $x\mapsto p.x$. Then the limit in $\dtop$ is the achronal
flow having as $0$-skeleton\index{$0$-skeleton} the set of $p$-adic
integers $\Z_p$ and the limit in $\cattopn$ is a totally disconnected
achronal $1$-category.

\bth\label{ss} The inclusion functor
$\widetilde{\omega}:\dtop\longrightarrow \cattopn$ has a right adjoint
that will be denoted by $\widetilde{D}$. In particular, this implies
that the canonical inclusion functor $\dtop\longrightarrow \cattopn $
preserves colimits. Moreover, one has $\widetilde{D}\circ
\widetilde{\omega}=\Id_{\dtop}$ and
\[\limproj_i X_i\iso \limproj_i\widetilde{D}\circ
\widetilde{\omega}\left(X_i\right)\iso \widetilde{D}
\left(\limproj_i \widetilde{\omega}\left(X_i\right)\right).\] \eth

If $\sets$ is the category of sets, then the forgetful functor
$\omega:\top\longrightarrow \sets$ has a left adjoint: the functor
$X\mapsto Dis\left(X\right)$ which maps a set $X$ to the discrete
space $Dis\left(X\right)$.  So
\[\top\left(Dis\left(X\right),Y\right)\iso \sets\left(X,\omega
\left(Y\right)\right).\]

\bpf Let $\C$ be an object of $\cattopn$. Then:
\begin{itemize}
\item Let $\widetilde{D}\left(\C\right)^0:=\C^0$ equipped with the discrete topology.
\item If $\left(\alpha,\beta\right)\in \widetilde{D}\left(\C\right)^0\p \widetilde{D}\left(\C\right)^0$,
let $\P_{\alpha,\beta}\widetilde{D}\left(\C\right)$ be the
subspace of $\P\C$ of execution paths $x$ such that $s(x)=\alpha$ and
$t(x)=\beta$ equipped with the Kelleyfication of the relative
topology.
\item Let $\P \widetilde{D}\left(\C\right)=\bigsqcup_{\left(\alpha,\beta\right)\in \widetilde{D}\left(\C\right)^0\p \widetilde{D}\left(\C\right)^0} \P_{\alpha,\beta}\widetilde{D}\left(\C\right)$ with an obvious definition of
the source map $s$, the target map $t$ and the composition law
$*$.
\end{itemize}
Let $f\in
\dtop\left(X,\widetilde{D}\left(Y\right)\right)$. Then the
composite $X^0\longrightarrow
\widetilde{D}\left(Y\right)^0\longrightarrow Y^0$ is continuous.  And
for any $\alpha,\beta\in X^0$, $\P_{\alpha,\beta}X\longrightarrow
\P_{f\left(\alpha\right),f\left(\beta\right)}Y\longrightarrow Y$ is
continuous as well. Reciprocally, a map $g\in
\cattopn\left(\widetilde{\omega}\left(X\right),Y\right)$ provides
$g^0\in
\top\left(\widetilde{\omega}\left(X\right)^0,Y^0\right)\iso
\set\left(\omega\circ\widetilde{\omega}\left(X\right)^0,\omega\left(Y^0\right)\right)$ since
$\widetilde{\omega}\left(X\right)^0$ is a discrete space and provides a continuous map {
\[\P g\in \top\left(\P\widetilde{\omega}\left(X\right),\P Y\right)\iso \top\left(\bigsqcup_{(\alpha,\beta)} \P_{\alpha,\beta}X,\P Y\right)\longrightarrow \prod_{(\alpha,\beta)}\top\left(\P_{\alpha,\beta}X,\P\widetilde{D} Y\right).\]}
Hence the natural bijection
\[\dtop\left(X,\widetilde{D}\left(Y\right)\right)\iso \cattopn\left(\widetilde{\omega}\left(X\right),Y\right).\]
\epf

\subsection{Tensor product of non-contracting topological $1$-ca\-te\-go\-ries}

The purpose of this section is the construction of a closed symmetric
monoidal structure on $\cattopn$.  Let \[\tcattopn(Y,Z)\] be the set
$\cattopn\left(Y,Z\right)\subset \ttop\left(Y,Z\right)$ equipped with
the Kelleyfication of the relative topology induced by that of
$\ttop\left(Y,Z\linebreak[0]\right)$.

\bp\label{tenseurn} Let $X$ and $Y$ be two objects of $\cattopn$.
There exists a unique structure of topological $1$-category
 $X\ot Y$ on the topological space  $\left(X^0\sqcup \P X\right)\p \left(Y^0\sqcup \P Y\right)$ such that
\begin{enumerate}
\item $\left(X\ot Y\right)^0=X^0\p Y^0$ .
\item $\P\left(X\ot Y\right)= \left(\P X\p \P X\right)\sqcup \left(X^0 \p \P Y\right) \sqcup \left(\P X\p Y^0\right)$.
\item $s\left(x,y\right)=\left(s(x),s(y)\right)$, $t\left(x,y\right)=\left(t(x),t(y)\right)$, $\left(x,y\right)*\left(z,t\right)=\left(x*z,y*t\right)$.
\end{enumerate}
 \ep

\bpf Obvious. \epf

\bp\label{composition} Let $X$ and $Y$ be two objects of
$\cattopn$. Let $f$ and $g$ be two morphisms in $\cattopn$ from
$\vI\ot X$ to $Y$. Let us suppose that for any $y\in Y$,
$f\left(1\ot y\right)=g\left(0\ot y\right)$. Then for any $y\in
X$, the following equality holds
\[ f\left([0,1]\ot s(y)\right)* g\left([0,1]\ot y\right)= f\left([0,1]\ot y\right)* g\left([0,1]\ot t(y)\right)\]
Denote the common value by $\left(f*g\right)\left([0,1]\ot
y\right)$. Let \[\left(f*g\right)\left(0\ot y\right)=f\left(0\ot
y\right)\] and \[\left(f*g\right)\left(1\ot y\right)=g\left(1\ot
y\right).\] Then $f*g$ yields an element of $\cattopn\left(\vI\ot
X,Y\right)$ and one has moreover
$\left(f*g\right)*h=f*\left(g*h\right)$. At last, this composition
yields a continuous map from the fiber product
\[\tcattopn\left(\vI\ot X,Y\right)\p_{\tcattopn\left(X,Y\right)} \tcattopn\left(\vI\ot X,Y\right)\]
given by the inclusions $\{0\}\subset \vI$ and $\{1\}\subset \vI$
to  $\tcattopn\left(\vI\ot X,Y\right)$. \ep

\bpf First of all, one has
\begin{alignat*}{2}
&f\left([0,1]\ot s(y)\right)* g\left([0,1]\ot y\right) \\
&= f\left([0,1]\ot s(y)\right)*g\left(0\ot y\right)*g\left([0,1]\ot t(y)\right)&&\ \hbox{since $g$ morphism of $\cattopn$}\\
&= f\left([0,1]\ot s(y)\right)*f\left(1\ot y\right)*g\left([0,1]\ot t(y)\right)&&\ \hbox{by hypothesis}\\
&= f\left([0,1]\ot y\right)*g\left([0,1]\ot t(y)\right)&&\
\hbox{since $f$ morphism of $\cattopn$}
\end{alignat*}
The equalities \[\left(f*g\right)\left(0\ot
x*y\right)=\left(f*g\right)\left(0\ot
x\right)*\left(f*g\right)\left(0\ot y\right)\] and
\[\left(f*g\right)\left(1\ot x*y\right)=\left(f*g\right)\left(1\ot
x\right)*\left(f*g\right)\left(1\ot y\right)\] are trivial.
Because of the symmetries, it remains to check that
\[\left(f*g\right)\left([0,1]\ot x*y\right)=\left(f*g\right)\left(0\ot x\right)*\left(f*g\right)\left([0,1]\ot y\right)\]
to get $f*g\in \cattopn\left(\vI\ot X,Y\right)$. And one has \beas
\left(f*g\right)\left([0,1]\ot x*y\right)&=& f\left([0,1]\ot \left(x*y\right)\right)*  g\left([0,1]\ot t\left(x*y\right)\right)\\
&=& f\left(0\ot x\right)*f\left([0,1]\ot y\right)* g\left(0\ot t(y)\right) * g\left([0,1]\ot t(y)\right)\\
&=& f\left(0\ot x\right)*f\left([0,1]\ot y\right)* f\left(1\ot t(y)\right) * g\left([0,1]\ot t(y)\right)\\
&=& f\left(0\ot x\right)*f\left([0,1]\ot y\right)*g\left([0,1]\ot t(y)\right)\\
&=& \left(f*g\right)\left(0\ot
x\right)*\left(f*g\right)\left([0,1]\ot y\right). \eeas At last,
one has to check that $\left(f*g\right)*h=f*\left(g*h\right)$.
Once again, the equalities
\[\left(\left(f*g\right)*h\right)\left(0\ot x\right)=\left(f*\left(g*h\right)\right)\left(0\ot x\right)\] and \[\left(\left(f*g\right)*h\right)\left(1\ot x\right)=\left(f*\left(g*h\right)\right)\left(1\ot x\right)\]
are trivial. And one has \beas
\left(\left(f*g\right)*h\right)\left([0,1]\ot x\right)&=& \left(f*g\right)\left([0,1]\ot s(x)\right)*h\left([0,1]\ot x\right)\\
&=& f\left([0,1]\ot s(x)\right)*g\left([0,1]\ot s(x)\right)*h\left([0,1]\ot x\right)\\
&=& f\left([0,1]\ot s(x)\right)*\left(g*h\right)\left([0,1]\ot x\right)\\
&=& \left(f*\left(g*h\right)\right)\left([0,1]\ot x\right). \eeas
The continuity of $*$ is due to the fact that we are working
exclusively with compactly generated topological spaces. \epf

\bth\label{na} The tensor product of $\cattopn$ is a closed
symmetric monoi\-dal structure, that is there exists a bifunctor
\[ \homcat:\cattopn\p \cattopn \longrightarrow \cattopn\]
contravariant with respect to the first argument and covariant
with respect to the second argument such that one has the natural
bijection of sets \[\cattopn\left(X\ot Y,Z\right)\iso
\cattopn\left(X,\homcat\left(Y,Z\right)\right)\] for any
topological $1$-categories $X$, $Y$ and $Z$. \eth

\bpf \

\fbox{
\begin{minipage}{.8\textwidth}
   Construction of $\homcat\left(Y,Z\right)$
\end{minipage}
}
\begin{enumerate}
\item $\homcat\left(Y,Z\right)^0:=\tcattopn\left(Y,Z\right)$
\item $\P\homcat\left(Y,Z\right):=\tcattopn\left(\vI\ot Y,Z\right)$
\item the source map and target map are induced respectively by
the morphisms $\{0\}\subset \vI$ and $\{1\}\subset \vI$
\item the composition law is defined by Proposition~\ref{composition}.
\end{enumerate}

\fbox{
\begin{minipage}{0.8\textwidth}
Construction of the set map $\Phi:\cattopn\left(X\ot
Y,Z\right)\longrightarrow
\cattopn\left(X,\homcat\left(Y,Z\right)\right)$ (with $f\in \cattopn\left(X\ot Y,Z\right)$)
\end{minipage}
}
\begin{enumerate}
\item for $x\in X^0$, $\Phi\left(f\right)\left(x\right)$ is the morphism of flows
from $Y$ to $Z$ defined by \begin{itemize}
\item $\Phi\left(f\right)\left(x\right)\left(y\right)=f\left(x\ot y\right)$. \end{itemize}
\item for $x\in \P X$, $\Phi\left(f\right)\left(x\right)$ is the morphism of flows from $\vI\ot Y$ to $Z$ defined
by \begin{itemize}
\item $\Phi\left(f\right)\left(x\right)\left(0\ot y\right)=f\left(s\left(x\right)\ot y\right)$
\item $\Phi\left(f\right)\left(x\right)\left(1\ot y\right)=f\left(t\left(x\right)\ot y\right)$
\item $\Phi\left(f\right)\left(x\right)\left([0,1]\ot y\right)=f\left(x\ot y\right)$.
\end{itemize}
\end{enumerate}

\fbox{
\begin{minipage}{0.8\textwidth}
Construction of the set map
$\Psi:\cattopn\left(X,\homcat\left(Y,Z\right)\right)\longrightarrow
\cattopn\left(X\ot Y,Z\right)$ (with $g\in
\cattopn\left(X,\homcat\left(Y,Z\right)\right)$)
\end{minipage}
}
\begin{enumerate}
\item $\Psi\left(g\right)\left(x_0\ot y\right)=g\left(x_0\right)\left(y\right)$ for $\left(x_0\ot y\right)\in X^0\p Y$
\item $\Psi\left(g\right)\left(x\ot y\right)=g\left(x\right)\left([0,1]\ot y\right)$
for $\left(x\ot y\right)\in \P X\p Y$.
\end{enumerate}

\fbox{
\begin{minipage}{0.8\textwidth}
$\Phi\left(f\right)\left(x *
x'\right)=\Phi\left(f\right)\left(x\right)*\Phi\left(f\right)\left(x'\right)$ (with $x,x'\in \P X$)
\end{minipage}
}
\begin{enumerate}
\item $\Phi\left(f\right)\left(x * x'\right)\left(0\ot y\right)=f\left(s\left(x\right)\ot y\right)$
\item $\Phi\left(f\right)\left(x * x'\right)\left([0,1]\ot y\right)=f\left((x * x')\ot y\right)$
\item $\Phi\left(f\right)\left(x * x'\right)\left(1\ot y\right)=f\left(t\left(x'\right)\ot y\right)$
\end{enumerate}

\fbox{
\begin{minipage}{0.8\textwidth}
$\Phi\left(f\right)\left(s\left(x\right)\right)=s\left(\Phi\left(f\right)\left(x\right)\right)$ and
$\Phi\left(f\right)\left(t\left(x\right)\right)=t\left(\Phi\left(f\right)\left(x\right)\right)$
\end{minipage}
}
\begin{enumerate}
\item $\Phi\left(f\right)\left(s\left(x\right)\right)=f\left(s\left(x\right)\ot -\right)= s\left(\Phi\left(f\right)\left(x\right)\right)$
\item  $\Phi\left(f\right)\left(t\left(x\right)\right)=f\left(t\left(x\right)\ot -\right)= t\left(\Phi\left(f\right)\left(x\right)\right)$.
\end{enumerate}

\fbox{
\begin{minipage}{0.8\textwidth}
$\Psi\left(g\right)\left(\left(x_0\ot y\right)*\left(x_0\ot
y'\right)\right)=\Psi\left(g\right)\left(x_0\ot
y\right)*\Psi\left(g\right)\left(x_0\ot y'\right)$ (with $g\in \dtop\left(X,\homcat\left(Y,Z\right)\right)$,
$x_0\in X^0$, $y,y'\in \P Y$)
\end{minipage}
}
\beas \Psi\left(g\right)\left(\left(x_0\ot y\right)*\left(x_0\ot
y'\right)\right)&=& \Psi\left(g\right)\left(\left(x_0\ot
(y*y')\right)\right)\\
&=&g\left(x_0\right)\left(y*y'\right)\\
&=&g\left(x_0\right)\left(y\right)*g\left(x_0\right)\left(y'\right)\\
&=&\Psi\left(g\right)\left(x_0\ot
y\right)*\Psi\left(g\right)\left(x_0\ot y'\right). \eeas

\fbox{
\begin{minipage}{0.8\textwidth}
$\Psi\left(g\right)\left(\left(x_0\ot y\right)*\left(x\ot y'\right)\right)=\Psi\left(g\right)\left(x_0\ot y\right)*\Psi\left(g\right)\left(x\ot y'\right)$ (with $g\in \dtop\left(X,\homcat\left(Y,Z\right)\right)$,
$x_0\in X^0$, $x\in \P X$, $y,y'\in \P Y$)
\end{minipage}
}
\beas
\Psi\left(g\right)\left(\left(x_0\ot y\right)*\left(x\ot y'\right)\right) &=& \Psi\left(g\right)\left(x\ot (y*y')\right) \\
&=& g\left(x\right)\left([0,1]\ot (y*y')\right) \\
&=&g\left(x\right)\left(0\ot y\right)*g\left(x\right)\left([0,1]\ot y'\right)\\
&=& \Psi\left(g\right)\left(x_0\ot y\right)*\Psi\left(g\right)\left(x\ot y'\right) \eeas

\fbox{
\begin{minipage}{0.8\textwidth}
$s\left(\Psi\left(g\right)\left(x\ot y\right)\right)=\Psi\left(g\right)\left(s\left(x\right)\ot s\left(y\right)\right)$ (with $x\in \P X$, $y\in Y$)
\end{minipage}
}
\beas s\left(\Psi\left(g\right)\left(x\ot y\right)\right)&=&
s\left(g\left(x\right)\left([0,1]\ot y\right)\right) \\ &=& g\left(x\right)\left(s\left([0,1]\ot y\right)\right) \\ &=& g\left(x\right)\left(0 \ot s\left(y\right)\right)
\\ &=& \left(s\left(g\left(x\right)\right)\right)\left(s\left(y\right)\right) \\ &=& \left(g\left(s\left(x\right)\right)\right)\left(s\left(y\right)\right) \\ &=& \Psi\left(g\right)\left(s\left(x\right)\ot s\left(y\right)\right) \eeas

\fbox{
\begin{minipage}{0.8\textwidth}
$s\left(\Psi\left(g\right)\left(x_0\ot y\right)\right)=\Psi\left(g\right)\left(x_0\ot s\left(y\right)\right)$ (with $x_0\in X^0$, $y\in Y$)
\end{minipage}
}
\beas s\left(\Psi\left(g\right)\left(x_0\ot y\right)\right)&=&
s\left(g\left(x_0\right)\left(y\right)\right) \\ &=& g\left(x_0\right)\left(s\left(y\right)\right)\\
&=& \Psi\left(g\right)\left(x_0\ot s\left(y\right)\right)
\eeas

\fbox{
\begin{minipage}{0.8\textwidth}
$\Phi\circ
\Psi=\Id_{\cattopn\left(X,\homcat\left(Y,Z\right)\right)}$
\end{minipage}
}

Let $x_0\in X^0$ and $y\in Y$. Then
\[\Phi\left(\Psi\left(g\right)\right)\left(x_0\right)\left(y\right)=\Psi\left(g\right)\left(x_0\ot y\right)=g\left(x_0\right)\left(y\right)\]
therefore
$\Phi\left(\Psi\left(g\right)\right)\left(x_0\right)=g\left(x_0\right)$.
And for $x\in \P X$,
\begin{enumerate}
\item $\Phi\left(\Psi\left(g\right)\right)\left(x\right)\left(0\ot y\right)=\Psi\left(g\right)\left(s(x)\ot y\right)=g\left(s(x)\right)\left(y\right)$
\item $\Phi\left(\Psi\left(g\right)\right)\left(x\right)\left(1\ot y\right)=\Psi\left(g\right)\left(t(x)\ot y\right)=g\left(t(x)\right)\left(y\right)$
\item $\Phi\left(\Psi\left(g\right)\right)\left(x\right)\left([0,1]\ot y\right)=\Psi\left(g\right)\left(x\ot y\right)=g\left(x\right)\left([0,1]\ot y\right)$.
\end{enumerate}

\fbox{
\begin{minipage}{0.8\textwidth}
$\Psi\circ
\Phi=\Id_{\cattopn\left(X\ot Y,Z\right)}$
\end{minipage}
}

With $f\in
\cattopn\left(X\ot Y,Z\right)$, $x_0\in X^0$ and $y\in Y$, one has
\[\Psi\left(\Phi\left(f\right)\right)\left(x_0\ot
y\right)=\Phi\left(f\right)\left(x_0\right)\left(y\right)=f\left(x_0\ot
y\right)\] and for $x\in \P X$,
\[\Psi\left(\Phi\left(f\right)\right)\left(x\ot y\right)=\Phi\left(f\right)\left(x\right)\left([0,1]\ot y\right)=f\left(x
\ot y\right).\]

\fbox{
\begin{minipage}{0.8\textwidth}
The continuity of $\Phi\left(f\right)$
\end{minipage}
}
\begin{enumerate}
\item The continuity of $\Phi\left(f\right)^0:X^0\longrightarrow \tcattopn\left(Y,Z\right)$ because
\[\Phi\left(f\right)^0\in \top\left(X^0,\ttop\left(Y,Z\right)\right)\iso \top\left(X^0\p Y,Z\right).\]
\item The continuity of $\P\Phi\left(f\right):\P X\longrightarrow \tcattopn\left(\vI\ot Y,Z\right)$
because \[\P\Phi\left(f\right)\in \top\left(\P X,\ttop\left(\vI\p
Y,Z\right)\right)\iso \top\left(\P X\p \vI\p Y,Z\right).\]
\end{enumerate}

\fbox{
\begin{minipage}{0.8\textwidth}
The continuity of $\Psi\left(g\right)$
\end{minipage}
}

The continuity of $\Psi\left(g\right)$ comes again from the
canonical bijections of sets
\[\top\left(X^0,\ttop\left(Y,Z\right)\right)\iso \top\left(X^0\p Y,Z\right)\]
and
\[\top\left(\P X,\ttop\left(\vI\p Y,Z\right)\right)\iso \top\left(\P X\p \vI\p Y,Z\right)\]
and also from the fact that the underlying topological space of a
given $1$-category $X$ is homeomorphic to the disjoint sum of
topological spaces $X^0\sqcup \P X$. This completes the proof. \epf

\begin{cor} Let $X$ and  $Y$  be two topological $1$-categories. Then one
has the homeomorphisms
\[\tcattopn(\liminj_i X_i,Y)\iso \limproj_i \tcattopn(X_i,Y)\]
and
\[\tcattopn(X,\limproj_i Y_i)\iso \limproj_i \tcattopn(X,Y_i)\]
for any colimit $\liminj_i X_i$ and any limit $\limproj_i Y_i$ in
$\cattopn$.
\end{cor}

In both following calculations, one uses the fact that the
following natural homeomorphism holds in $\cattopn$:
$\left(\limproj_i X_i\right)^0\iso \limproj_i
\left(X_i^0\right)$. The latter homeomorphism may be false in
$\dtop$ since the $0$-skeleton\index{$0$-skeleton} is always discrete
in the latter category.

\bpf
One has: \beas
\tcattopn\left(\liminj_i X_i,Y\right)&\iso & \left(\homcat\left(\liminj_i X_i,Y\right)\right)^0\\&\iso & \left(\limproj_i \homcat\left(X_i,Y\right)\right)^0\\
&\iso & \limproj_i \left(\homcat\left(X_i,Y\right)\right)^0\\
&\iso & \limproj_i \tcattopn\left( X_i,Y\right) \eeas and \beas
\tcattopn\left( X,\limproj_iY_i\right)&\iso & \left(\homcat\left(X,\limproj_i Y_i\right)\right)^0\\&\iso & \left(\limproj_i \homcat\left(X,Y_i\right)\right)^0\\
&\iso & \limproj_i \left(\homcat\left(X,Y_i\right)\right)^0\\
&\iso & \limproj_i \tcattopn\left(X,Y_i\right). \eeas
 \epf

\subsection{Important consequence for the category of flows}

As an application of the preceding results, one proves the
following crucial theorem:

\bth\label{commute} Let $\tdtop\left(X,Y\right)$ be the set of
morphisms of flows from $X$ to $Y$ equipped with the
Kelleyfication of the compact-open\index{compact-open topology}
topology. Then the mapping \[\left(X,Y\right)\mapsto
\tdtop\left(X,Y\right)\] induces  a functor from $\dtop\p \dtop$
to $\top$ which is contravariant with respect to $X$ and covariant
with respect to $Y$. Moreover:
\begin{enumerate}
\item One has  the natural homeomorphism \[\tdtop\left(\liminj_i X_i,Y\right)\iso
\limproj_i \tdtop\left(X_i,Y\right)\] for any colimit $\liminj_i
X_i$ in $\dtop$.
\item One has the natural homeomorphism \[\tdtop\left(X,\limproj_i Y_i\right)\iso
\limproj_i\tdtop\left(X, Y_i\right)\] for any finite limit
$\limproj_i X_i$ in $\dtop$.
\end{enumerate}
\eth

The functor $\tdtop(X,-)$ cannot commute with any limit. Indeed, with
$X=\{0\}$, one has $\tdtop(X,Y)\iso Y^0$ as space. However, a limit of
a diagram of discrete topological space may be totally disconnected
without being discrete.

This is the reason why we make the distinction between the
\textit{set} of morphisms $\dtop(X,Y)$ from a flow $X$ to a flow $Y$
and the \textit{space} of morphisms $\tdtop(X,Y)$ from a flow $X$ to a
flow $Y$.

\bpf Since $\widetilde{\omega}$ preserves colimits by
Theorem~\ref{ss}, one has: \beas
\tdtop\left(\liminj_i X_i,Y\right) &\iso & \tcattopn\left(\widetilde{\omega}\left(\liminj_i X_i\right),\widetilde{\omega}\left(Y\right)\right)\\
&\iso & \tcattopn\left(\liminj_i\widetilde{\omega}\left( X_i\right),\widetilde{\omega}\left(Y\right)\right)\\
&\iso &\limproj_i\tcattopn\left(\widetilde{\omega}\left( X_i\right),\widetilde{\omega}\left(Y\right)\right)\\
&\iso &\limproj_i\tdtop\left(X_i,Y\right) \eeas Since
$\widetilde{\omega}$ preserves finite limits by Theorem~\ref{rr},
one has: \beas
\tdtop\left(X,\limproj_i Y_i\right) &\iso & \tcattopn\left(\widetilde{\omega}\left(X\right),\widetilde{\omega}\left(\limproj_i Y_i\right)\right)\\
 &\iso & \tcattopn\left(\widetilde{\omega}\left(X\right),\limproj_i \widetilde{\omega}\left(Y_i\right)\right)\\
&\iso &  \limproj_i \tcattopn\left(\widetilde{\omega}\left(X\right),\widetilde{\omega}\left(Y_i\right)\right)\\
&\iso &  \limproj_i \tdtop\left(X,Y_i\right). \eeas \epf

One does not need actually the previous machinery of tensor
product of $1$-categories to prove the isomorphism of topological
spaces
\[\tdtop\left(X,\limproj_i Y_i\right)\iso
\limproj_i\tdtop\left(X, Y_i\right)\] for any finite limit
$\limproj_i Y_i$ of $\dtop$. Indeed one sees that the forgetful
functor $X\mapsto X^0\sqcup \P X$ from $\dtop$ to $\top$ induces the
inclusion of topological spaces {\small
\[\tdtop\left(X,\limproj_i Y_i\right)\subset \ttop\left(X^0,\limproj_i  Y_i^0\right)\p \prod_{\left(\alpha,\beta\right)\in X^0\p X^0}\ttop\left(\P_{\alpha,\beta} X,\limproj_i  \P_{\alpha_i,\beta_i} Y_i\right)\]}
where $\alpha_i$ (resp. $\beta_i$) is the image of $\alpha$
(resp. $\beta$) by the composite $X^0 \longrightarrow \limproj_i
Y_i^0\longrightarrow Y_i^0$. Since the right member of the above
inclusion is isomorphic to {\[\limproj_i\left(
\ttop\left(X^0,Y_i^0\right)\p \prod_{\left(\alpha,\beta\right)\in
X^0\p X^0}\ttop\left(\P_{\alpha,\beta} X, \P_{\alpha_i,\beta_i}
Y_i\right) \right)
\]}
then the conclusion follows.

On the contrary, the forgetful functor $X\mapsto X^0\sqcup \P X$ from
$\dtop$ to $\top$ does not commute at all with colimits, even the
finite ones, because colimits in $1$-categories may create execution
paths. So the tensor product of $1$-categories seems to be required to
establish the other homeomorphism.

\section{Flow as a canonical colimit of globes and
points}\label{canonique}

In the sequel, one will implicitely use the category
$\diag\left(\dtop\right)$ of diagrams of flows. The
objects are the functor $D:I\longrightarrow \dtop$ where $I$ is a
small category.  A morphism from a diagram $D:I\longrightarrow
\dtop$ to a diagram $E:J\longrightarrow \dtop$ is a functor
$\phi:I\longrightarrow J$ together with a natural transformation
$\mu:D\longrightarrow E\circ \phi$. A morphism of diagram
$\left(\phi,\mu\right):D\longrightarrow E$ gives rise to a morphism
of flows $\liminj D\longrightarrow \liminj E$. Since $\dtop$ is
complete and cocomplete, then $\diag\left(\dtop\right)$ is
complete and cocomplete as well \cite{MR57:9788}.

In this section, we prove that any flow is the colimit in a
canonical way of globes and points. This technical tool will be
used in the sequel of the paper.

\bth \label{point-globe} Any flow is the colimit in $\dtop$ of
points and globes in a canonical way, i.e. there exists for any
flow $X$ a diagram $\mathbb{D}\left(X\right)$ of flows containing
only points, globes and concatenations of globes such that the
mapping $X\mapsto \mathbb{D}\left(X\right)$ is functorial and
such that $X\iso \liminj \mathbb{D}\left(X\right)$ in a canonical
way. \eth

\bpf Let $X$ be a flow and let $\alpha$, $\beta$ and $\gamma$ be
three points (not necessarily distinct) of its
$0$-skeleton\index{$0$-skeleton}. Consider the diagram of
Figure~\ref{UFO} where the map
\[\glob\left(\P_{\alpha,\beta}X\p \P_{\beta,\gamma}X\right)\longrightarrow
\glob\left(\P_{\alpha,\beta}X\right)*\glob\left(\P_{\beta,\gamma}X\right)\]
is induced by the map $(x,y)\mapsto x*y$ (where $*$ is the free
concatenation) and where the map \[\glob\left(\P_{\alpha,\beta}X\p
\P_{\beta,\gamma}X\right) \longrightarrow
\glob\left(\P_{\alpha,\gamma}X\right)\] is induced by the
composition law of $X$. Then consider the diagram
$\mathbb{D}\left(X\right)$ obtained by concatening all diagrams as
that of Figure~\ref{UFO}. It is constructed as follows:
\begin{itemize}
\item the underlying small category $\mathcal{I}\left(X\right)$ of $\mathbb{D}\left(X\right)$
is the free category generated by the set of objects $X^0\cup
X^0\p X^0\cup X^0\p X^0\p X^0\p \{0,1\}$ and by the arrows
\begin{itemize}
\item $i_1^{\alpha,\beta}:\alpha\longrightarrow \left(\alpha,\beta\right)$ and
$i_2^{\alpha,\beta}:\beta\longrightarrow \left(\alpha,\beta\right)$
\item $r^{\alpha,\beta,\gamma}:\left(\alpha,\beta,\gamma,0\right)\longrightarrow
\left(\alpha,\beta,\gamma,1\right)$ and
$p^{\alpha,\beta,\gamma}:\left(\alpha,\beta,\gamma,0\right)\longrightarrow
\left(\alpha,\gamma\right)$
\item $j_1^{\alpha,\beta,\gamma}:\alpha\longrightarrow
\left(\alpha,\beta,\gamma,0\right)$ and
$j_3^{\alpha,\beta,\gamma}:\gamma\longrightarrow
\left(\alpha,\beta,\gamma,0\right)$
\item $k_1^{\alpha,\beta,\gamma}:\alpha\longrightarrow
\left(\alpha,\beta,\gamma,1\right)$,
$k_2^{\alpha,\beta,\gamma}:\beta\longrightarrow
\left(\alpha,\beta,\gamma,1\right)$ and
$k_3^{\alpha,\beta,\gamma}:\gamma\longrightarrow
\left(\alpha,\beta,\gamma,1\right)$
\item $h_1^{\alpha,\beta,\gamma}:\left(\alpha,\beta\right)\longrightarrow
\left(\alpha,\beta,\gamma,1\right)$ and
$h_3^{\alpha,\beta,\gamma}:\left(\beta,\gamma\right)\longrightarrow
\left(\alpha,\beta,\gamma,1\right)$
\end{itemize}
\item $\mathbb{D}\left(X\right)\left(\alpha\right)=\{\alpha\}$,
$\mathbb{D}\left(X\right)\left(\alpha,\beta\right)=\glob\left(\P_{\alpha,\beta}X\right)$
\item
$\mathbb{D}\left(X\right)\left(\alpha,\beta,\gamma,0\right)=\glob\left(\P_{\alpha\beta}X\p
\P_{\beta,\gamma}X\right)$
\item
$\mathbb{D}\left(X\right)\left(\alpha,\beta,\gamma,1\right)=\glob\left(\P_{\alpha\beta}X\right)*\glob\left(\P_{\beta,\gamma}X\right)$
\item $\mathbb{D}\left(X\right)\left(i_1^{\alpha,\beta}\right)$ is the canonical inclusion
\[\{\alpha\}\longrightarrow \glob\left(\P_{\alpha,\beta}X\right)\]
\item $\mathbb{D}\left(X\right)\left(i_2^{\alpha,\beta}\right)$ is the canonical inclusion
\[\{\beta\}\longrightarrow \glob\left(\P_{\alpha,\beta}X\right)\]
\item $\mathbb{D}\left(X\right)\left(r^{\alpha,\beta,\gamma}\right)$ is the canonical
projection \[\glob\left(\P_{\alpha,\beta}X\p
\P_{\beta,\gamma}X\right)\longrightarrow
\glob\left(\P_{\alpha,\beta}X\right)*\glob\left(\P_{\beta,\gamma}X\right)\]
sending $\left(x,y\right)$ to $x*y$.
\item $\mathbb{D}\left(X\right)\left(p^{\alpha,\beta,\gamma}\right):
\glob\left(\P_{\alpha\beta}X\p
\P_{\beta,\gamma}X\right)\longrightarrow
\glob\left(\P_{\alpha,\gamma}X\right)$ is the morphism  induced by
the composition law of $X$
\item $\mathbb{D}\left(X\right)\left(j_1^{\alpha,\beta,\gamma}\right)$ (resp.
$\mathbb{D}\left(X\right)\left(j_3^{\alpha,\beta,\gamma}\right)$)
is the canonical inclusion from $\{\alpha\}$ (resp. $\{\gamma\}$)
to
\[\glob\left(\P_{\alpha,\beta}X\p \P_{\beta,\gamma}X\right)\]
\item $\mathbb{D}\left(X\right)\left(k_1^{\alpha,\beta,\gamma}\right)$ (resp.
$\mathbb{D}\left(X\right)\left(k_2^{\alpha,\beta,\gamma}\right)$,
$\mathbb{D}\left(X\right)\left(k_3^{\alpha,\beta,\gamma}\right)$)
is the canonical inclusion from $\{\alpha\}$ (resp. $\{\beta\}$,
$\{\gamma\}$) to
\[\glob\left(\P_{\alpha\beta}X\right)*\glob\left(\P_{\beta,\gamma}X\right)\]
\item $\mathbb{D}\left(X\right)\left(h_1^{\alpha,\beta,\gamma}\right)$ is the canonical
inclusion \[\glob\left(\P_{\alpha,\beta}X\right)\longrightarrow
\glob\left(\P_{\alpha,\beta}X\right)*\glob\left(\P_{\beta,\gamma}X\right)\]
\item $\mathbb{D}\left(X\right)\left(h_3^{\alpha,\beta,\gamma}\right)$ is the canonical
inclusion \[\glob\left(\P_{\beta,\gamma}X\right)\longrightarrow
\glob\left(\P_{\alpha,\beta}X\right)*\glob\left(\P_{\beta,\gamma}X\right)\]
\end{itemize}
Let $T$ be a flow. Let $f:\liminj_{i\in\mathcal{I}\left(X\right)}
\mathbb{D}\left(X\right)\left(i\right)\longrightarrow T$ be a
morphism of flows. Notice that all morphisms of flows in the
diagram $\mathbb{D}\left(X\right)$ are source and target
preserving. So $f$ yields a well-defined set map $g^0$ from
$Y^0=X^0$ to $T^0$. Moreover the morphism $f$ yields a continuous
map $f_{\alpha,\beta}:\P_{\alpha,\beta}X\longrightarrow
\P_{f\left(\alpha\right),f\left(\beta\right)}T$ and for any
$\left(\alpha,\beta,\gamma\right)\in Y^0\p Y^0\p Y^0$, a
continuous map \[f_{\alpha,\beta,\gamma,0}:\P_{\alpha,\beta}X\p
\P_{\beta,\gamma}X\longrightarrow
\P_{f\left(\alpha\right),f\left(\gamma\right)}T\] and another
continuous map \[f_{\alpha,\beta,\gamma,1}:\P_{\alpha,\beta}X\cup
\P_{\beta,\gamma}X\cup \P_{\alpha,\beta}X\p
\P_{\beta,\gamma}X\longrightarrow
\P_{f\left(\alpha\right),f\left(\beta\right)}T \cup
\P_{f\left(\beta\right),f\left(\gamma\right)}T\cup
\P_{f\left(\alpha\right),f\left(\gamma\right)}T\] which satisfy
various commutativity conditions. In particular all these maps
define a unique continuous map
$g_{\alpha,\beta}:\P_{\alpha,\beta}X\longrightarrow
\P_{f\left(\alpha\right),f\left(\beta\right)}T$ thanks to
$h_1^{\alpha,\beta,\gamma}$ and $h_3^{\alpha,\beta,\gamma}$
(these latter being inclusions). For $x\in \P_{\alpha,\beta}X$
and $y\in \P_{\beta,\gamma}X$, one has: \beas
g_{\alpha,\gamma}\left(x*y\right) &=& g_{\alpha,\gamma}\left(
p^{\alpha,\beta,\gamma}\left(x,y\right)\right)\\
&=& f_{\alpha,\beta,\gamma,0}\left(x,y\right)\\
&=& f_{\alpha,\beta,\gamma,1}\left(r^{\alpha,\beta,\gamma}\left(x,y\right)\right)\\
&=& f_{\alpha,\beta,\gamma,1}\left(x*y\right)\\
&=& f_{\alpha,\beta,1}\left(x\right)*f_{\beta,\gamma,1}\left(y\right) \hbox{ since $f$ morphism of flows !}\\
&=&
f_{\alpha,\beta,\gamma,1}\left(h_1^{\alpha,\beta,\gamma}\left(x\right)\right)*f_{\alpha,\beta,\gamma,1}
\left(h_3^{\alpha,\beta,\gamma}\left(y\right)\right)\\
&=&
f_{\alpha,\beta}\left(x\right)*f_{\beta,\gamma}\left(y\right)\\
&=&
g_{\alpha,\beta}\left(x\right)*g_{\beta,\gamma}\left(y\right)\eeas
So $g$ yields a well-defined morphism of flows from $X$ to $T$.
Conversely from a morphism of flows from $X$ to $T$, one can
construct a morphism of flows from
$\liminj_{i\in\mathcal{I}\left(X\right)}
\mathbb{D}\left(X\right)\left(i\right)$ to $T$. So one has the
natural bijection of sets
\[\dtop\left(\liminj_{i\in\mathcal{I}\left(X\right)} \mathbb{D}\left(X\right)\left(i\right),T\right)\iso \dtop\left(X,T\right)\]
Hence by Yoneda\index{Yoneda's lemma}, the flow $X$ is the
colimit of this diagram and moreover everything is canonical.
The  functoriality of $\mathbb{D}$ is obvious. \epf

{
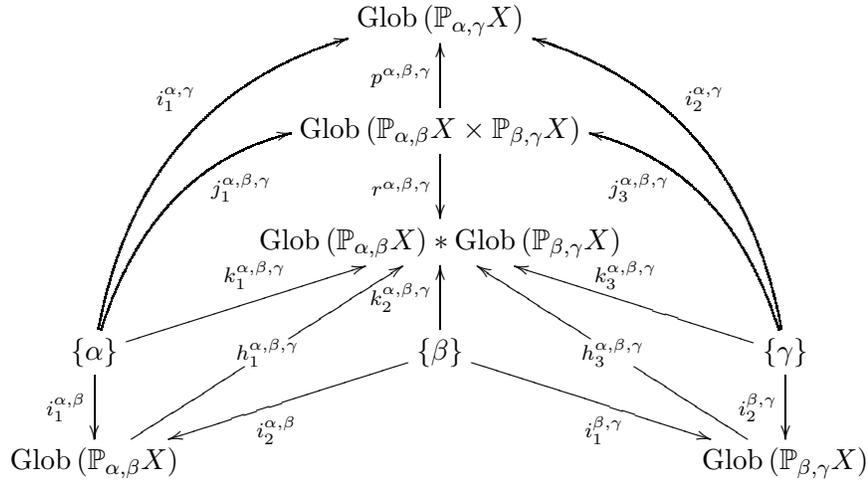
\begin{figure}
\[
\xymatrix{
& \glob\left(\P_{\alpha,\gamma}X\right) & \\
& \glob\left(\P_{\alpha,\beta}X\p \P_{\beta,\gamma}X\right)\fu{p^{\alpha,\beta,\gamma}} \fd{r^{\alpha,\beta,\gamma}}&\\
& \glob\left(\P_{\alpha,\beta}X\right)*\glob\left(\P_{\beta,\gamma}X\right)& \\
{\{\alpha\}}\fd{i_1^{\alpha,\beta}}\ar@{->}[ru]^{k_1^{\alpha,\beta,\gamma}}\ar@/^35pt/[ruu]_{j_1^{\alpha,\beta,\gamma}}\ar@/^35pt/[ruuu]^{i_1^{\alpha,\gamma}}
& \fu{k_2^{\alpha,\beta,\gamma}}
{\{\beta\}}\ar@{->}[ld]^{i_2^{\alpha,\beta}}\ar@{->}[rd]_{i_1^{\beta,\gamma}}
&
\ar@{->}[lu]_{k_3^{\alpha,\beta,\gamma}}{\{\gamma\}\fd{i_2^{\beta,\gamma}}\ar@/_35pt/[luu]^{j_3^{\alpha,\beta,\gamma}}\ar@/_35pt/[luuu]_{i_2^{\alpha,\gamma}}}\\
\glob\left(\P_{\alpha,\beta}X\right)\ar@{->}[ruu]|{h_1^{\alpha,\beta,\gamma}}
& &
\glob\left(\P_{\beta,\gamma}X\right)\ar@{->}[luu]|{h_3^{\alpha,\beta,\gamma}}
}
\]
\caption{The flow $X$ as a colimit of globes and points}
\label{UFO}
\end{figure}}

\begin{cor} \label{reduction}
Let $P\left(X\right)$ be a statement depending on a flow $X$ and
satisfying the following property: if $D:\mathcal{I}\longrightarrow
\dtop$ is a diagram of flows such that for any object $i$ of
$\mathcal{I}$, $P\left(D\left(i\right)\right)$ holds, then
$P\left(\liminj D\right)$ holds. Then the following assertions are
equivalent:
\begin{enumerate}
\item[(i)] The statement $P\left(X\right)$ holds for any flow $X$ of $\dtop$.
\item[(ii)] The statements $P\left(\{*\}\right)$ and $P\left(\glob\left(Z\right)\right)$ hold for any
object $Z$ of $\top$.
\end{enumerate}
\end{cor}

\bpf The implication $\left(i\right)\Longrightarrow
\left(ii\right)$ is obvious. Conversely if $\left(ii\right)$
holds, then
\[P\left(\glob\left(Z_1\right)*\glob\left(Z_2\right)\right)\] holds
for any topological spaces $Z_1$ and $Z_2$ since
$\glob\left(Z_1\right)*\glob\left(Z_2\right)$ is the colimit of
the diagram of flows
\[\xymatrix{\{*\}\fr{*\mapsto 1}\fd{*\mapsto 0} & \glob\left(Z_1\right)\\ \glob\left(Z_2\right) & }\]
containing only points and globes. The proof is complete with
Theorem~\ref{point-globe}. \epf

\section{S-homotopy in $\dtop$} \label{sectionS}

\subsection{Synchronized morphism of flows}

\bd A morphism of flows $f:X\longrightarrow Y$ is said {\rm
synchronized} if and only if it induces a bijection of sets
between the $0$-skeleton\index{$0$-skeleton} of $X$ and the
$0$-skeleton\index{$0$-skeleton} of $Y$. \ed

\subsection{S-homotopy of flows}
\label{sectionSS}

We mimick here the definition of the S-homotopy\index{S-homotopy}
relation for globular complexes \cite{diCW}.

\bd \label{defS} Let $f$ and $g$ be two morphisms of flows from
$X$ to $Y$. Then $f$ and $g$ are {\rm S-homotopic} or {\rm S-homotopy
equivalent} if there exists a continuous map $H:[0,1]\p
X\longrightarrow Y$ such that, with $H\left(u,-\right)=H_u$, for any
$u\in [0,1]$, $H_u$ is a morphism of flows from $X$ to $Y$ with
$H_0=f$ and $H_1=g$. In particular, this implies that $f$ and $g$
coincide on the $0$-skeleton\index{$0$-skeleton} $X^0$ of $X$ and that
for any $x_0\in X^0$, for any $u\in [0,1]$,
$f\left(x_0\right)=H\left(u,x_0\right)=g\left(x_0\right)$. This
situation is denoted by $f\sim_{S} g$. This defines an equivalence
relation on the set $\dtop\left(X,Y\right)$. \ed

Following Proposition~\ref{can}, one then obtains the natural
definition

\bd Two elements of the path space $\P X$ of a flow $X$ are said
{\rm S-homotopic} if the corresponding morphisms of flows from $\vI$
to $X$ are S-homotopy\index{S-homotopy} equivalent. \ed

\bd Two flows $X$ and $Y$ are {\rm S-homotopic} or {\rm S-homotopy equivalent} if there exists
two morphisms $f:X\longrightarrow Y$ and $g:Y\longrightarrow X$ such
that $f\circ g \sim_{S} \Id_Y$ and $g\circ f \sim_{S} \Id_X$. The maps
$f$ and $g$ are called (reciprocal) S-homotopy\index{S-homotopy}
equivalences. The S-homotopy\index{S-homotopy} relation is obviously
an equivalence relation.  We say that the mapping $g$ is a {\rm
S-homotopic inverse} of $f$. \ed

Because of the discreteness of the $0$-skeleton\index{$0$-skeleton} of
any flow, a S-homotopy\index{S-homotopy} equivalence is necessarily
synchronized.

\bp \label{caracflow} Let $f$ and $g$ be two morphisms of flows
from $X$ to $Y$. Then $f$ and $g$ are S-homotopic if and only if there
exists a continuous map \[h\in
\top\left([0,1],\tdtop\left(X,Y\right)\right)\] such that
$h\left(0\right)=f$ and $h\left(1\right)=g$. \ep

\bpf Let $H:[0,1]\p X \longrightarrow Y$ be the
S-homotopy\index{S-homotopy} from $f$ to $g$. Then $H$ provides an
element of $\top\left([0,1],\ttop\left(X,Y\right)\right)$ which is by
definition of a S-homotopy\index{S-homotopy} also an element of
$\top\left([0,1],\tdtop\left(X,Y\right)\right)$.  Conversely, an
element $h$ of
\[\top\left([0,1],\tdtop\left(X,Y\right)\right)\] yields an element
of
$\top\left([0,1],\ttop\left(X,Y\right)\right)\iso
\top\left([0,1]\p X,Y\right)$ which is by construction a
S-homotopy\index{S-homotopy} from $f$ to $g$. \epf

\subsection{Pairing $\boxtimes$ between a topological space and a flow}

\begin{nota} Let $U$ be a topological space. Let $X$ be a flow. The flow
$\{U,X\}_S$ is defined as follows:
\begin{enumerate}
\item The $0$-skeleton of $\{U,X\}_S$ is $X^0$.
\item For $\alpha,\beta\in X^0$, the topological space $\P_{\alpha,\beta}\{U,X\}_S$
is $\ttop(U,\P_{\alpha,\beta}X)$.
\item For $\alpha,\beta,\gamma\in X^0$, the composition law
\[*:\P_{\alpha,\beta}\{U,X\}_S\p
\P_{\beta,\gamma}\{U,X\}_S\longrightarrow\P_{\alpha,\gamma}\{U,X\}_S\]
is the composite {\small
\[\P_{\alpha,\beta}\{U,X\}_S\p \P_{\beta,\gamma}\{U,X\}_S\iso
\ttop\left(U,\P_{\alpha,\beta}X\p
\P_{\beta,\gamma}X\right)\longrightarrow
\ttop\left(U,\P_{\alpha,\gamma}X\right)\]} induced by the
composition law of $X$.
\end{enumerate}
\end{nota}

If $U=\varnothing$ is the empty set, then $\{\varnothing,Y\}_S$
is the flow having the same $0$-skeleton as $Y$ and exactly one
non-constant execution path between two points of $Y^0$.

\bth Let $U$ be a topological space. The mapping $Y\mapsto
\{U,Y\}_S$ yields a functor from $\dtop$ to itself. Moreover one
has
\begin{enumerate}
\item one has the natural isomorphism of flows $\{U,\limproj_i
X_i\}_S\iso \limproj_i \{U,X_i\}_S$
\item if $Y=Y^0$, then $\{U,Y\}_S=Y$
\item if $U$ and $V$ are two topological spaces, then
$\{U\p V,Y\}_S\iso \{U,\{V,Y\}_S\}_S$.
\end{enumerate}
\eth

\bpf   The functoriality of $\{U,-\}_S$ is obvious.  Following
the proof of Theorem~\ref{lim-colim}, it is clear that the
functor $\{U,-\}_S$ does preserve limits in $\dtop$. By
definition, \[\{U\p V,Y\}_{S}^0\iso
\{U,\{V,Y\}_S\}_{S}^0\iso Y^0\] and for $\alpha,\beta\in
Y^0$, one has \[\P_{\alpha,\beta}\{U\p V,Y\}_S=\ttop\left(U\p
V,\P_{\alpha,\beta}Y\right)\] and
\[\P_{\alpha,\beta}\{U,\{V,Y\}_S\}_S=\ttop\left(U,\ttop\left(V,\P_{\alpha,\beta}Y\right)\right).\]
Therefore $\{U\p V,Y\}_S\iso \{U,\{V,Y\}_S\}_S$. \epf

\bth \label{existenceflow} Let $U$ be a topological space. The
functor $\{U,-\}_S$ has a left adjoint which will be denoted by
$U\boxtimes -$. Moreover:
\begin{enumerate}
\item one has the natural isomorphism of flows \[U\boxtimes \left(\liminj_i X_i\right)
\iso \liminj_i \left(U\boxtimes X_i\right)\]
\item there is a natural isomorphism of flows $\{*\}\boxtimes Y\iso Y$
\item if $Z$ is a topological space, one has the natural
isomorphism of flows \[U\boxtimes \glob\left(Z\right)\iso
\glob\left(U\p Z\right)\]
\item for any flow $X$ and any topological space $U$, one has
the natural bijection of sets \[\left(U\boxtimes X\right)^0\iso
X^0\]
\item if $U$ and $V$ are two topological spaces, then $\left(U\p V\right)\boxtimes Y\iso
U\boxtimes \left(V\boxtimes Y\right)$ as flows
\item for any flow $X$, $\varnothing \boxtimes X\iso X^0$.
\end{enumerate}
\eth

If $u\in U$, the image of $x\in X$ by the canonical morphism of flows
$X\longrightarrow \{u\}\boxtimes X\longrightarrow U\boxtimes X$ is
denoted by $u\boxtimes x$.

\bpf In the category of $\dtop$, let us start with the class of solutions
$f:Z\longrightarrow \{U,Y\}_S$ for $f$ running over the set
$\dtop\left(Z,\{U,Y\}_S\right)$ and for $Y$ running over the
class of  flows. Consider the commutative diagram
\[\xymatrix{Z\fr{f} \ar@{->}[rd]^-{f}& \{U,Y'\}_S \fd{} \\ & \{U,Y\}_S}\]
where $Y'$ is the subflow generated by the elements of
$f\left(Z\right)\left(U\right)$ and where the vertical map is induced
by the inclusion $Y'\subset Y$.  So one still has a set of solutions
by considering only the flows $Y$ such that the cardinal $\card(Y)$ of
the underlying set satisfies $\card(Y)\leq
\aleph_0\p \card(Z)\p \card(U)$. Let $\{Z_i,i\in I\}$ be the set of
isomorphism classes of flows whose underlying set is of cardinal less
than $\aleph_0\p \card(Z)\p \card(U)$. Then $\card(I)\leq
2^{\left(\aleph_0\p \card(Z)\p \card(U)\right)^5}$ so $I$ is a set.
Then the class $\bigcup_{i\in I}\dtop\left(Z,\{U,Z_i\}_S\right)$ is a
set as well and one gets a set of solutions.  The first assertion is
then clear using Theorem~\ref{ssc}. One has $\dtop\left(\{*\}\boxtimes
X,Y\right)\iso
\dtop\left(X,\{\{*\},Y\}_S\right)\iso \dtop\left(X,Y\right)$ so
by Yoneda\index{Yoneda's lemma} $\{*\}\boxtimes X\iso X$ for any
flow $X$. One has \beas
\dtop\left(U\boxtimes \glob\left(Z\right),Y\right)&\iso & \dtop\left(\glob\left(Z\right),\{U,Y\}_S\right)\\
&\iso & \bigsqcup_{\left(\alpha,\beta\right)\in Y^0\p Y^0} \top\left(Z,\ttop\left(U,\P_{\alpha,\beta}Y\right)\right)\\
&\iso & \bigsqcup_{\left(\alpha,\beta\right)\in Y^0\p Y^0} \top\left(U\p Z, \P_{\alpha,\beta}Y\right)\\
&\iso & \dtop\left(\glob\left(U\p Z\right),Y\right) \eeas So by
Yoneda\index{Yoneda's lemma} $U\boxtimes \glob\left(Z\right)\iso
\glob\left(U\p Z\right)$. One has \beas
\dtop\left(U\boxtimes \{*\},Y\right)&\iso & \dtop\left( \{*\},\{\{*\},Y\}_S\right)\\
&\iso & \{\{*\},Y\}_{S}^0\\
&\iso & \dtop\left(\{*\},Y\right) \eeas so by
Yoneda\index{Yoneda's lemma}, $U\boxtimes \{*\}\iso \{*\}$. Hence
$\left(U\boxtimes X\right)^0\iso X^0$ if $X$ is a point or a
globe. Hence the result by Corollary~\ref{reduction}. One has \beas
\dtop\left(\left(U\p V\right)\boxtimes X,Y\right)&\iso & \dtop\left(X,\{V\p U,Y\}_S\right)\\
&\iso & \dtop\left(X,\{V,\{U,Y\}_S\}_S\right)\\
&\iso & \dtop\left(V\boxtimes X,\{U,Y\}_S\right)\\
&\iso & \dtop\left(U\boxtimes\left(V\boxtimes X\right),Y\right)
\eeas so by Yoneda\index{Yoneda's lemma} $\left(U\p
V\right)\boxtimes X\iso U\boxtimes\left(V\boxtimes X\right)$. \epf

Take a flow $X$ and a topological space $U$. One knows that $X$ is the
colimit in a canonical way of points and globes
(Theorem~\ref{point-globe}). Since $U\boxtimes \{*\}\iso \{*\}$ and
$U\boxtimes \glob\left(Z\right)\iso \glob\left(U\p Z\right)$, and
since the functor $U\boxtimes -$ commutes with colimits, one can
represent $U\boxtimes X$ as the colimit of the diagram of
Figure~\ref{UFO2} with an obvious definition of the arrows (in
particular $r^{\alpha,\beta,\gamma}_U$ uses the diagonal
$U\longrightarrow U\p U$).

\begin{bigcenter}
{
\begin{figure}
{\footnotesize
\[
\xymatrix{
& \glob\left(U\p\P_{\alpha,\gamma}X\right) & \\
& \glob\left(U\p\P_{\alpha,\beta}X\p \P_{\beta,\gamma}X\right)\fu{p^{\alpha,\beta,\gamma}_U} \fd{r^{\alpha,\beta,\gamma}_U}&\\
& \glob\left(U\p\P_{\alpha,\beta}X\right)*\glob\left(U\p\P_{\beta,\gamma}X\right)& \\
{\{\alpha\}}\fd{i_{1,U}^{\alpha,\beta}}\ar@{->}[ru]^{k_{1,U}^{\alpha,\beta,\gamma}}\ar@/^35pt/[ruu]_{j_{1,U}^{\alpha,\beta,\gamma}}\ar@/^35pt/[ruuu]^{i_{1,U}^{\alpha,\gamma}}
& \fu{k_{2,U}^{\alpha,\beta,\gamma}}
{\{\beta\}}\ar@{->}[ld]^{i_{2,U}^{\alpha,\beta}}\ar@{->}[rd]_{i_{1,U}^{\beta,\gamma}}
&
\ar@{->}[lu]_{k_{3,U}^{\alpha,\beta,\gamma}}{\{\gamma\}\fd{i_{2,U}^{\beta,\gamma}}\ar@/_35pt/[luu]^{j_{3,U}^{\alpha,\beta,\gamma}}\ar@/_35pt/[luuu]_{i_{2,U}^{\alpha,\gamma}}}\\
\glob\left(U\p\P_{\alpha,\beta}X\right)\ar@{->}[ruu]|{h_{1,U}^{\alpha,\beta,\gamma}}
& &
\glob\left(U\p\P_{\beta,\gamma}X\right)\ar@{->}[luu]|{h_{3,U}^{\alpha,\beta,\gamma}}
}
\]
}
\caption{Representation of the flow $U\boxtimes X$} \label{UFO2}
\end{figure}
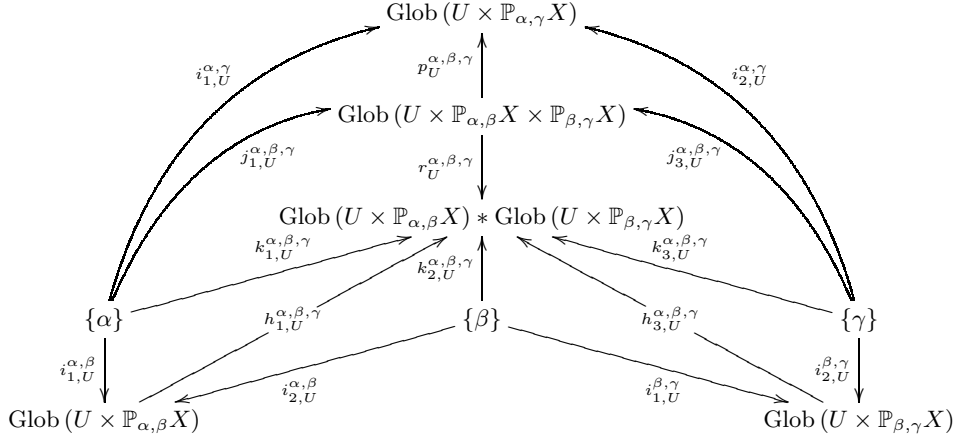}
\end{bigcenter}

\subsection{Cylinder functor for the S-homotopy of flows}

\bth\label{finpreuve} Let $U$ be a connected non-empty
topological space. Let $X$ and $Y$ be two flows. Then one has a
natural bijection of sets
\[\dtop\left(X,\{U,Y\}_S\right)\iso \top\left(U,\tdtop\left(X,Y\right)\right)\] and so
\[\dtop\left(U\boxtimes X,Y\right)\iso \top\left(U,\tdtop\left(X,Y\right)\right).\]
\eth

\bpf It suffices to prove the first bijection by Theorem~\ref{existenceflow}.
Let \[f\in \dtop\left(X,\{U,Y\}_S\right).\] Then $f$ induces a set map
from $X^0$ to $Y^0$ (but $Y^0\iso \top(U,Y^0)$ since $U$ is connected)
and a continuous map from $\P X$ to $\ttop(U,\P Y)$. So one has the
inclusion of sets \[i_1:\dtop\left(X,\{U,Y\}_S\right)\subset
\top\left(X,\ttop\left(U,Y\right)\right).\] The inclusion of
sets $\tdtop\left(X,Y\right)\longrightarrow\ttop\left(X,Y\right)$
induces an inclusion of sets \[i_2:\top\left(U,\tdtop\left(X,Y\right)\right)\longrightarrow
\top\left(U,\ttop\left(X,Y\right)\right).\] But
$\top\left(X,\ttop\left(U,Y\right)\right)\iso
\top\left(U,\ttop\left(X,Y\right)\right)$.  And it is then easy
to see that $i_1$ and $i_2$ have the same
image. So the sets $\dtop\left(X,\{U,Y\}_S\right)$ and
\[\top\left(U,\tdtop\left(X,Y\right)\right)\] are bijective. \epf

\bd Let $\C$ be a category. A {\rm cylinder} is a functor
$I:\C\longrightarrow \C$ together with natural transformations
$i_0,i_1:\id_\C\longrightarrow I$ and $p:I\longrightarrow \id_\C$ such
that $p\circ i_0$ and $p\circ i_1$ are the identity natural
transformation. \ed

\begin{cor}[Cylinder functor] \label{Scyl}
The mapping $X\mapsto [0,1]\boxtimes X$ induces a functor from $\dtop$
to itself which is a cylinder functor\index{cylinder functor} with the
natural transformations $e_i:\{i\}\boxtimes -
\longrightarrow [0,1]\boxtimes -$ induced by the inclusion maps
$\{i\}\subset [0,1]$ for $i\in\{0,1\}$ and with the natural
transformation $p:[0,1]\boxtimes -\longrightarrow \{0\}\boxtimes -$
induced by the constant map $[0,1]\longrightarrow \{0\}$. Moreover,
two morphisms of flows $f$ and $g$ from $X$ to $Y$ are S-homotopic if
and only if there exists a morphism of flows $H:[0,1]\boxtimes
X\longrightarrow Y$ such that $H\circ e_0=f$ and $H\circ
e_1=g$. Moreover $e_0\circ H\sim_{S} \Id$ and $e_1\circ H\sim_{S}
\Id$.
\end{cor}

\bpf Consequence of Theorem~\ref{finpreuve}, Proposition~\ref{caracflow} and
Theorem~\ref{existenceflow} and of the connectedness of $[0,1]$. \epf

\section{Explicit description of $U\boxtimes X$}\label{boxex}

\bp \label{disparition} Let $X$ and $Y$ be two flows. Let $U$ be
a topological space. Then one has a bijection between the elements of
$\dtop\left(U\boxtimes X,Y\right)$ and the elements $f$ of
$\sets\left(X^0,Y^0\right)\p \top\left(U\p \P X,Y\right)$ such that
\begin{itemize}
\item $f\left(X^0\right)\subset Y^0$
\item $f\left(U\p \P X\right)\subset \P Y$
\item for any $u\in U$, $f\left(u,x*y\right)=f\left(u,x\right)*f\left(u,y\right)$ if $x,y\in \P X$ and if $tx=sy$
\item for any $u\in U$, $s(f\left(u,x\right))=f\left(s(x)\right)$ and $t(f\left(u,x\right))=f\left(t(x)\right)$ if $x\in \P X$.
\end{itemize}
\ep

\bpf The set $\dtop\left(U\boxtimes X,Y\right)$ is isomorphic to the
set $\dtop\left(X,\{U,Y\}_S\right)$, hence the result. \epf

\bp\label{freeflow} The forgetful functor $\upsilon$ from $\dtop$
to the category of diagrams $D$ of topological spaces over the small
category $\xymatrix{1\fr{s} & 0 & \fl{t} 1}$ such that $D(0)$ is a
discrete topological space has a left adjoint called the free flow
generated by the diagram. \ep

\bpf The forgetful functor preserves limits because of the
construction of the limit in $\dtop$. Let $D$ be a diagram over
$\xymatrix{1\fr{s} & 0 & \fl{t} 1}$ with $D\left(0\right)$
discrete. Let us start from the class of solutions $f:D\longrightarrow
\upsilon\left(X\right)$ when $X$ runs over the class of flows and for
a given $X$ where $f$ runs over
$\dtop\left(D,\upsilon\left(X\right)\right)$. Then one can replace $X$
by the subflow generated by the finite composition of elements of
$f\left(D\right)$. So one can suppose that the cardinal $\card(X)$ of
$X$ satisfies $\card(X)\leq \aleph_0\p \card(D)$ where $\card(D)$ is
the cardinal of $D$. By choosing one equivalence class of flows for
the class of flows $X$ such that $\card(X)\leq \aleph_0\p
\card(D)$, one has obtained a set of solutions. Hence the result by Theorem~\ref{ssc}. \epf

\begin{cor}\label{disparition1}
Let $U$ be a topological space. Let $X$ be a flow. Then the flow
$U\boxtimes X$ is the free flow generated by the diagram $D$ of spaces
defined by $D\left(1\right)=U\p \P X$, $D\left(0\right)=X^0$,
$s\left(u,x\right)=s(x)$, $t\left(u,x\right)=t(x)$ divided by the
identifications
$\left(u,x\right)*\left(u,y\right)=\left(u,x*y\right)$. \end{cor}

\bpf The identifications generates an equivalence with closed graph since the
composition law of $X$ is continuous. Therefore the quotient equipped
with the final topology is still weak Hausdorff, and therefore
compactly generated. This is then a consequence of Yoneda's
lemma. \epf

\section{S-homotopy extension property}\label{SDEP}

\bd Let $i:A\longrightarrow X$ be a synchronized morphism of flows and let $Y$ be
a flow. The morphism $i:A\longrightarrow X$ satisfies {\rm the
S-homotopy extension property}\index{S-homotopy extension property}
for $Y$ if for any morphism $f:X\longrightarrow Y$ and any
S-homotopy\index{S-homotopy} $h:[0,1]\boxtimes A\longrightarrow Y$
such that for any $a\in A$, $h\left(0\boxtimes
a\right)=f\left(i\left(a\right)\right)$, there exists a
S-homotopy\index{S-homotopy} $H:[0,1]\boxtimes X\longrightarrow Y$
such that for any $x\in X$, $H\left(0\boxtimes
x\right)=f\left(x\right)$ and for any $\left(t,a\right)\in [0,1]\p A$,
$H\left(t\boxtimes i\left(a\right)\right)=h\left(t\boxtimes
a\right)$. \ed

\bd A synchronized morphism $i:A\longrightarrow X$ satisfies the
{\rm S-homotopy extension property} if $i:A\longrightarrow X$
satisfies the S-homotopy\index{S-homotopy} extension
property\index{S-homotopy extension property} for any flow $Y$.
\ed

Let $i:A\longrightarrow X$ be a morphism of flows and let
$\left(Y,B\right)$ be a pair of topological spaces. Then one can
consider the pushout
\[
\xymatrix{B\boxtimes A \fr{} \fd{\Id_B\boxtimes i}& Y\boxtimes A \fd{}\\
B\boxtimes X \fr{} & \cocartesien \left(Y\boxtimes
A\right)\sqcup_{B\boxtimes A} \left(B\boxtimes X\right) }
\]
The commutativity of the diagram
\[
\xymatrix{B\boxtimes A \fr{} \fd{\Id_B\boxtimes i}& Y\boxtimes A \fd{}\\
B\boxtimes X \fr{} &  Y\boxtimes X }
\]
provides a canonical morphism of flows $\left(Y,B\right)\boxtimes
i: \left(Y\boxtimes A\right)\sqcup_{B\boxtimes A}
\left(B\boxtimes X\right)\longrightarrow Y\boxtimes X$ such that
\begin{itemize}
\item[] $\left(\left(Y,B\right)\boxtimes i\right)\left(y\boxtimes a\right)=y\boxtimes i\left(a\right)$,
\item[] $\left(\left(Y,B\right)\boxtimes i\right)\left(b\boxtimes a\right)=b\boxtimes i\left(a\right)$,
\item[] $\left(\left(Y,B\right)\boxtimes i\right)\left(b\boxtimes x\right)=b\boxtimes x$
\end{itemize}
with $x\in X$, $y\in Y$, $a\in A$ and $b\in B$.

\begin{nota} The
morphism of flows $\left([0,1],\{0\}\right)\boxtimes i$ will be
denoted by $\psi(i)$. \end{nota}

\bth\label{retract} Let $i:A\longrightarrow X$ be a morphism of
flows. Then the following assertions are equivalent:
\begin{enumerate}
\item the morphism $i$ satisfies the S-homotopy extension property
\item the morphism of flows $\psi(i)$ has a retract $r$,
that is to say there exists a morphism of flows
\[r:[0,1]\boxtimes X\longrightarrow \left([0,1]\boxtimes A\right)\sqcup_{\{0\}\boxtimes A} \left(\{0\}\boxtimes X\right)\]
such that $r\circ \psi(i)=\Id_{\left([0,1]\boxtimes
A\right)\sqcup_{\{0\}\boxtimes A} \left(\{0\}\boxtimes X\right)}$.
\end{enumerate}
\eth

\bpf Giving two morphisms of flows $f:X\longrightarrow Y$ and
$h:[0,1]\boxtimes A\longrightarrow Y$ such that $h\left(0\boxtimes
a\right)=f\left(i\left(a\right)\right)$ for any $a\in A$ is equivalent
to giving a morphism of flows still denoted by $h$ from
$\left([0,1]\boxtimes A\right)\sqcup_{\{0\}\boxtimes A}
\left(\{0\}\boxtimes X\right)$ to $Y$. The
S-homotopy\index{S-homotopy} extension problem for $i$ has then always
a solution if and only for any morphism of flows
$h:\left([0,1]\boxtimes A\right)\sqcup_{\{0\}\boxtimes A}
\left(\{0\}\boxtimes X\right)\longrightarrow Y$, there exists a
morphism of flows $H:[0,1]\boxtimes X\longrightarrow Y$ such that
$H\circ \psi(i)=h$. Take $Y=\left([0,1]\boxtimes
A\right)\sqcup_{\{0\}\boxtimes A} \left(\{0\}\boxtimes X\right)$
and let $h$ be the identity map of $Y$. This yields the retract
$r$.  Conversely, let $r$ be a retract of $i$. Then $H:=h\circ r$
is always a solution of the S-homotopy\index{S-homotopy}
extension problem.  \epf

\bth\label{ndr} Let $(Z,\de Z)$ be a NDR pair of topological
spaces. Then the canonical morphism of flows $\glob(\de
Z)\longrightarrow \glob(Z)$ satisfies the S-homotopy extension property. \eth

\bpf Since $(Z,\de Z)$ is a NDR pair, then $[0,1]\p \de Z
\sqcup_{\{0\}\p \de Z} Z\longrightarrow  [0,1]\p Z$ has a retract.
Therefore the morphism of flows
\[\glob\left([0,1]\p \de Z \sqcup_{\{0\}\p \de Z} Z\right)\longrightarrow \glob\left([0,1]\p Z\right)\]
has a retract. But
\[\glob\left([0,1]\p \de Z \sqcup_{\{0\}\p \de Z} Z\right)\iso [0,1]\boxtimes \glob(\de Z)
\sqcup_{\{0\}\boxtimes \glob(\de Z)} \glob(Z)\] and
$\glob\left([0,1]\p Z\right)\iso [0,1]\boxtimes \glob(Z)$. The
proof is complete thanks to Theorem~\ref{retract}. \epf

\bth\label{limtopt} Let $U$ be a connected non empty space.
Let $X$ and $Y$ be two flows. Then there exists
a natural homeomorphism
\[\ttop(U,\tdtop(X,Y))\iso \tdtop(U\boxtimes X,Y).\]
\eth

\bpf We already know by Theorem~\ref{finpreuve} that there exists
a natural bijection
\[\top(U,\tdtop(X,Y))\iso \gltop(U\boxtimes X,Y).\]
Using the construction of $\boxtimes$, Corollary~\ref{reduction} and
Theorem~\ref{commute}, it suffices to prove the homeomorphism for
$X=X^0$ and $X=\glob(Z)$.  The space $\tdtop(X^0,Y)$ is the discrete
space of set maps $\set(X^0,Y^0)$ from $X^0$ to $Y^0$. Since $U$ is
connected, then $\ttop(U,\tdtop(X^0,Y))\iso
\set(X^0,Y^0)$. In the other hand,
$\tdtop(U\boxtimes X^0,Y)\iso \tdtop(X^0,Y)\iso \set(X^0,Y^0)$, hence
the result for $X^0$. At last, for any topological space $W$,
\beas
&& \top\left(W,\ttop\left(U,\tdtop\left(\glob\left(Z\right),Y\right)\right)\right)\\
&&\iso \top\left(W\p U,\tdtop\left(\glob\left(Z\right),Y\right)\right)\\
&&\iso \gltop\left(\left(W\p U\right)\boxtimes \glob\left(Z\right),Y\right)\\
&&\iso \gltop\left(\glob\left(W\p U\p Z\right),Y\right) \eeas and
$\top(W,\tdtop(U\boxtimes \glob(Z),Y)\iso \top(W,\tdtop(\glob(U\p
Z),Y))$. It is then easy to see that both \[\gltop(\glob(W\p U\p
Z),Y)\] and \[\top(W,\tdtop(\glob(U\p Z),Y))\] can be identified
to the same subset of $\top([0,1]\p W\p U\p Z,Y)$. Hence the
result by Yoneda. \epf

\bth\label{car} A morphism of flows $i:A\longrightarrow X$ satisfies
the S-homotopy extension property if and only if for any flow $Y$, the
continuous map $i^*:\tdtop(X,Y)\longrightarrow\tdtop(A,Y)$ is a
Hurewicz fibration. \eth

\bpf For any topological space $M$, one has
\[\top([0,1]\p M,\tdtop(A,Y))\iso \top (M,\ttop([0,1],\tdtop(A,Y)))\]
since $\top$ is cartesian closed and \[\top
(M,\ttop([0,1],\tdtop(A,Y)))\iso \top(M,\tdtop([0,1]\boxtimes A,Y))\] by
Theorem~\ref{limtopt}. Considering a  commutative diagram like
\[
\xymatrix{
\{0\}\p M\ar@{^{(}->}[d] \fr{\phi} & \tdtop(X,Y)\fd{i^*}\\
[0,1]\p M \fr{\psi}\ar@{-->}[ru]^-{k} & \tdtop(A,Y)}\]
is then equivalent to considering a commutative diagram of topological spaces
\[\xymatrix{
M \fd{}\fr{} & \tdtop(\{0\}\boxtimes X,Y) \fd{}\\
\tdtop([0,1]\boxtimes A,Y)\fr{} & \tdtop(\{0\}\boxtimes A,Y)}
\]
Using again Theorem~\ref{commute}, considering such
a commutative diagram is equivalent to considering a continuous map
$M\longrightarrow \tdtop(Mi,Y)$. Finding a continuous map $k$ making
both triangles commutative is equivalent to finding a commutative
diagram of the form
\[\xymatrix{
M \fr{\overline{\phi}} \fd{=}&  \tdtop(Mi,Y) \\
M \ar@{-->}[r]^-{\ell}& \tdtop([0,1]\boxtimes X,Y)\fu{\psi(i)^*}}
\]
If $i:A\longrightarrow X$ satisfies the S-homotopy extension property,
then $\psi(i):Mi\longrightarrow [0,1]\boxtimes X$ has a retract
$r:[0,1]\boxtimes X\longrightarrow Mi$.  Then take
$\ell=\overline{\phi}\circ r$. Conversely, if $\ell$ exists for any
$M$ and any $Y$, take $M=\{0\}$ and $Y=Mi$ and
$\overline{\phi}(0)=\id_{Mi}$. Then $\ell(0)$ is a retract of
$\psi(i)$. Therefore $i:A\longrightarrow X$ satisfies the S-homotopy
extension property. \epf

\begin{cor} \label{pushout-Scof} 
Let $i:A\longrightarrow X$ satisfy the S-homotopy extension property.
Let $f:A\longrightarrow Y$ be a morphism of flows.  Consider the
pushout in $\dtop$
\[\xymatrix{A \ar@{->}[r]^-{i} \ar@{->}[d]_-{f} & X \fd{} \\ Y \fr{j} & Z  \cocartesien}
\]
Then the canonical morphism from $Y$ to $Z$ satisfies the S-homotopy
extension property.  In other terms, the pushout of a morphism of
flows satisfying the S-homotopy extension property still satisfies the
S-homotopy extension property. \end{cor}

\begin{cor}  Let $Z$ be a compact space and let $\de
Z\subset Z$ be a compact subspace such that the canonical
inclusion is a NDR pair.  Let $U$ be a flow. Then the canonical
restriction map \[
\tdtop\left(\glob\left(Z\right),U\right)\longrightarrow
\tdtop\left(\glob\left(\de Z\right),U\right)\] is a Hurewicz fibration. \end{cor}

\section{Morphisms of flows inducing a closed inclusion of path spaces}\label{secincl}

\bp \cite{Ref_wH} If $j:X\longrightarrow Y$ and
$r:Y\longrightarrow X$ are two continuous maps with $r\circ j=\id$,
then $j$ is a closed inclusion and $r$ is a quotient map.
\ep

\begin{nota} Denote by $\mathcal{INC}$ the class of morphisms of flows
$f:X\longrightarrow Y$ such that $\P f:\P X\longrightarrow \P Y$ is a
closed inclusion of topological spaces.
\end{nota}

The purpose of this section is to collect some important examples of
morphisms of $\mathcal{INC}$. This section provides the necessary
preparatory lemmas for the use of the ``Small Object Argument''
further in this paper.

\bp Let $A$ be a flow. Then the morphism of flows
$\theta:A\longrightarrow [0,1]\boxtimes A$ defined by
$\theta(a)=1\boxtimes a$ belongs to $\mathcal{INC}$.  \ep

\bpf The mapping $U\mapsto U\boxtimes X$ for a given flow $X$ is
functorial with respect to $U$. So one can consider $r:[0,1]\boxtimes
A\longrightarrow A$ defined by $r(t\boxtimes a)=a$. \epf

\bd If $i:A\longrightarrow X$ is a morphism of flows, then the
{\rm mapping cylinder} $Mi$ of $i$ is defined by the pushout of flows
\[\xymatrix{
A \fr{a\mapsto 0\boxtimes a} \fd{i} & [0,1]\boxtimes A\fd{}\\
X \fr{} & \cocartesien Mi}
\]
\ed

\bp\label{theta} Let $i:A\longrightarrow X$ be a synchronized
morphism of flows. Then the canonical morphism of flows
$\theta:A\longrightarrow Mi$ such that $\theta(a)=1\boxtimes a$
belongs to $\mathcal{INC}$.  \ep

\bpf First of all, since $i$ is synchronized, one can consider
that $A^0=X^0$. Let \[Ni=([0,1]\boxtimes A) \sqcup_{A^0} X.\] Then
there exists a canonical morphism of flows $\phi:Ni\longrightarrow Mi$
which is constant on the $0$-skeleton and such that $\phi:\P
Ni\longrightarrow \P Mi$ is onto. Let us consider the equivalence
relation $\mathcal{R}$ on $\P Ni$ associated to $\phi$,
i.e. $x\mathcal{R}y$ if and only if $\phi(x)=\phi(y)$. The graph of
$\mathcal{R}$ is the inverse image of the diagonal of $\P Mi$. The
latter is closed in $\P Mi\p \P Mi$ since $\P Mi$ is a $k$-space which
is weak Hausdorff.  Therefore the graph of $\mathcal{R}$ is closed in
$\P Ni\p \P Ni$. So the quotient $\P Ni/\mathcal{R}$ equipped with the
final topology is still weak Hausdorff and thus a compactly generated
topological space. There exists a canonical continuous map $\P
Ni/\mathcal{R}\longrightarrow \P Mi$ which is an isomorphism of
sets. The topological space $\P Ni/\mathcal{R}$ yields a flow $Y$ and
a commutative diagram of flows
\[\xymatrix{
A \fr{a\mapsto 0\boxtimes a} \fd{i} & [0,1]\boxtimes A\fd{}\\
X \fr{} & Y}
\]
Therefore there exists a morphism of flows $Mi\longrightarrow Y$
because of the universal property satisfied by $Mi$.  So the bijection
$\sigma:\P Ni/\mathcal{R}\longrightarrow \P Mi$ is actually an
homeomorphism and the continuous map $\phi:\P Ni\longrightarrow \P Mi$
is a quotient map. The map $\theta:A\longrightarrow Mi$ is equal to
the composite
\[A\iso \{1\}\boxtimes A \longrightarrow Ni \longrightarrow Mi\]
One has $\P Ni\iso \P([0,1]\boxtimes A)\sqcup Z$ for some topological
space $Z$: the topological space $Z$ consists of all free compositions
of executions paths of $Ni$ containing an element of $X$. Therefore
the morphism of flows $A\longrightarrow Ni$ is a closed inclusion of
topological spaces. Let $g:Z\longrightarrow A$ be a set map such that
$\theta\circ g:Z\longrightarrow \P Mi$ is continuous. Let $F$ be a
closed subspace of $A$. Then $F$ is mapped to a closed subspace $G$ of
$\P Ni$. Since $G=\sigma^{-1}(\sigma(G))$, then $\sigma(G)$ is a
closed subspace of $\P Mi$. Therefore $g^{-1}(F)=(\theta\circ
g)^{-1}(\sigma(G))$ is a closed subspace of $Z$. Therefore $g$ is
continuous. \epf

\bth\label{exincl}
Let $i:A\longrightarrow X$ satisfy the S-homotopy extension property.
Then $i\in \mathcal{INC}$. \eth

\bpf We follow the proof of the fact that any Hurewicz cofibration
of compactly topological spaces is a closed inclusion given in the
appendix of \cite{Ref_wH}.

Let us consider the commutative diagram of flows
\[
\xymatrix{ A \fr{\theta} \fd{i} & Mi\fd{} \\
X \fr{i_1} & [0,1]\boxtimes X}
\]
where $\theta(a)=1\boxtimes a$ and $i_1(x)=1\boxtimes x$.  Then $i_1$
has a retract and therefore is a closed inclusion. The map $\theta$ is
a closed inclusion as well by Proposition~\ref{theta}.  Since $j$ has
a retract by Theorem~\ref{retract}, then $j\circ\theta$ is a closed
inclusion. moreover $i_1$ is one-to-one. Therefore $i$ is a closed
inclusion. \epf

\section{Smallness argument}\label{small}

Any ordinal can be viewed as a small category whose objects are the
elements of $\lambda$, that is the ordinal $\gamma<\lambda$, and where
there exists a morphism $\gamma\longrightarrow \gamma'$ if and only if
$\gamma\leq \gamma'$.

\bd Let $\C$ be a cocomplete category. Let $\lambda$ be an
ordinal. A {\rm $\lambda$-sequence} in $\C$ is a colimit-preserving
functor $X:\lambda\longrightarrow \C$. Since $X$ preserves colimits,
for all limit ordinals $\gamma<\lambda$, the induces map
$\liminj_{\beta<\gamma}X_\beta\longrightarrow X_\gamma$ is an
isomorphism. The morphism $X_0\longrightarrow\liminj X$ is called the
{\rm transfinite composition} of the $X_\gamma\longrightarrow
X_{\gamma+1}$. \ed

\bd Let $\kappa$ be a cardinal. An ordinal $\lambda$ is {\rm $\kappa$-filtered}
if for any $A\subset \lambda$ with $|A|\leq \kappa$ where $|A|$ is the
cardinal of $A$, then $\sup A<\lambda$.
\ed

\bd Let $\C$ be a cocomplete category. Let $\D$ be a collection
of morphisms of $\C$. Let $\kappa$ be a cardinal.  An object $A$ of
$\C$ is \textit{$\kappa$-small} with respect to $\D$ if for any
$\lambda$-sequence $X$ where $\lambda$ is a $\kappa$-filtered ordinal,
and where each arrow $X_\beta\longrightarrow X_{\beta+1}$ lies in $\D$
for $\beta<\lambda$, then one has the bijection
$\liminj_{\beta<\lambda} \C(A,X_\beta)\longrightarrow
\C(A,\liminj_{\beta<\lambda} X_\beta)$. We say that $A$ is
$\kappa$-small relative to $\mathcal{D}$ if it is $\kappa$-small
relative to $\D$ for some cardinal $\kappa$. \ed

\bd Let $\C$ be a cocomplete category. Let $I$ be a set of
morphisms of $\C$. Then a \textit{relative $I$-cell complex}
$f:A\longrightarrow B$ is a transfinite composition of pushouts of
elements of $I$. In other terms, there exists an ordinal
$\lambda$ and a $\lambda$-sequence $X:\lambda\longrightarrow \C$ such
that $f$ is the composition of $X$ and such that for each $\beta$
with $\beta+1<\lambda$, there is a pushout square as follows
\[
\xymatrix{ C_\beta\fr{} \fd{g_\beta} & X_\beta\fd{}\\ D_\beta
\fr{} & \cocartesien X_{\beta+1}}
\]
such that $g_\beta\in I$. We denote the collection of relative
$I$-cell complexes by $I$-cell. If $\varnothing$ is the initial object
of $\C$ and if $X$ is an object of $\C$ such that
$\varnothing\longrightarrow X$ is a relative $I$-cell complex, then one
says that $X$ is a $I$-cell complex. \ed

\bp\label{compact} Any flow $A$ is
$\sup(\aleph_0,\card(A))$-small relative to $\mathcal{INC}$ where
$\card(A)$ is the cardinal of the underlying topological space of $A$.
\ep

\bpf  One has a canonical one-to-one set map
\[\liminj_{\beta<\lambda} \dtop(A,X_\beta)\longrightarrow
\dtop(A,\liminj_{\beta<\lambda} X_\beta).\]
Let $f\in \dtop(A,\liminj_{\beta<\lambda} X_\beta)$. Since the
$0$-skeleton of a colimit of flows is the colimit of the
$0$-skeletons, then for any $a\in A^0$, $f(a)\in X_{\beta_a}^0$ for
some $\beta_a<\lambda$. There exists a canonical continuous map
$\liminj_{\beta<\lambda} \P X_\beta\longrightarrow
\P\left(\liminj_{\beta<\lambda} X_\beta\right)$ where
$\liminj_{\beta<\lambda} X_\beta$ is the colimit of the flows
$X_\beta$. Any element of $\P\left(\liminj_{\beta<\lambda}
X_\beta\right)$ is a finite composite $x_1*\dots*x_r$ of elements
$x_1\in X_{\beta_1},\dots,x_r\in X_{\beta_r}$ for some finite integer
$r$. Since $\lambda$ is $\sup(\aleph_0,\card(A))$-filtered, it is
$\aleph_0$-filtered. So $\beta=\sup(\beta_1,\dots,\beta_r)<\lambda$
and $x_1,\dots,x_r\in \P X_\beta$. So $x_1*\dots*x_r\in \P X_\beta$.
Therefore any execution path $x\in \P\left(\liminj_{\beta<\lambda}
X_\beta\right)$ belongs to some $\P X_{\beta_x}$ for some
$\beta_x<\lambda$.  Since $\lambda$ is
$\sup(\aleph_0,\card(A))$-filtered, it is
$\card(A)$-filtered. Therefore
$\sup(\beta_a,\dots,\beta_x)<\lambda$. So $f$ factors through a map
$g:A\longrightarrow X_\beta$ with $\beta<\lambda$. The map
$g:A\longrightarrow X_\beta$ is automatically continuous because all
continuous maps between path spaces are inclusions of topological
spaces.
\epf

\bd A morphism of flows $f:X\longrightarrow Y$ is a {\rm weak
S-homotopy equivalence} if $f$ is synchronized and if $f$ induces a
weak homotopy equivalence from $\P X$ to $\P Y$. \ed

\begin{nota} Let $\mathcal{S}$ be the subcategory of weak S-homotopy
equivalences. Let $I^{gl}$ be the set of morphisms of flows
$\glob(\mathbf{S}^{n-1})\longrightarrow \glob(\mathbf{D}^n)$ for $n\geq 0$. Let
$J^{gl}$ be the set of morphisms of flows
$\glob(\mathbf{D}^{n})\longrightarrow \glob([0,1]\p \mathbf{D}^n)$. Notice that all
arrows of  $\mathcal{S}$, $I^{gl}$ and $J^{gl}$ are
synchronized. At last, denote by $I^{gl}_+$ be the union of
$I^{gl}$ with the two morphisms of flows $R:\{0,1\}\longrightarrow
\{0\}$ and $C:\varnothing\subset \{0\}$. \end{nota}

\bp The domains of $I^{gl}_+$ are small relative to $I^{gl}_+$-cell.
 The domains of $J^{gl}$ are small relative to $J^{gl}$-cell. \ep

\bpf The inclusion maps $\mathbf{S}^{n-1}\subset \mathbf{D}^n$ and
$\mathbf{D}^{n}\subset [0,1]\p \mathbf{D}^n$ are NDR pairs. So any
pushout of a morphism of $I^{gl}\cup J^{gl}$ satisfies the
S-homotopy extension property by Corollary~\ref{pushout-Scof} and
Theorem~\ref{ndr}, and therefore is an element of $\mathcal{INC}$
by Theorem~\ref{exincl}. A pushout of
$C:\varnothing\longrightarrow \{0\}$ does not change the path
space. Therefore such a pushout is necessarily in
$\mathcal{INC}$. It remains to examine the case of a pushout of
$R:\{0,1\}\longrightarrow \{0\}$. Let us consider the pushout of
flows
\[\xymatrix{
\{0,1\}\fr{\phi}\fd{R} & X\fd{}\\
\{0\} \fr{} & Y \cocartesien}
\]
If $\phi(0)=\phi(1)$, then $\P X=\P Y$ and so there is nothing to
prove. Otherwise, if $\phi(0)\neq\phi(1)$, then $\P Y\iso \P X\sqcup
(\P_{.,\phi(1)}X\p \P_{\phi(0),.}X)\sqcup (\P_{.,\phi(1)}X\p
\P_{\phi(0),\phi(1)}X\p \P_{\phi(0),.}X)\sqcup
\dots$. Hence the conclusion by Proposition~\ref{compact}. \epf

\section{Reminder about model category}\label{remindermodel}

Some useful references for the notion of model\index{model
category} category are \cite{MR99h:55031} \cite{MR2001d:55012}. See also
\cite{ref_model1} \cite{ref_model2}.

If $\C$ is a category, one denotes by $Map(\C)$ the category whose
objects are the morphisms of $\C$ and whose morphisms are the
commutative squares of $\C$.

In a category $\C$, an object $x$ is \textit{a retract} of an
object $y$ if there exists $f:x\longrightarrow y$ and $g:y\longrightarrow
x$ of $\C$ such that $g\circ f=\Id_x$. A \textit{functorial
factorization} $(\alpha,\beta)$ of $\C$ is a pair of functors
from $Map(\C)$ to $Map(\C)$ such that for any $f$ object of
$Map(\C)$, $f=\beta(f)\circ \alpha(f)$.

\bd Let $i:A\longrightarrow B$ and $p:X\longrightarrow Y$ be maps in a
category $\C$. Then $i$ has the \textit{left lifting property}
(LLP) with respect to $p$ (or $p$ has the \textit{right lifting
property} (RLP) with respect to $i$) if for any commutative square
\[
\xymatrix{
A\fd{i} \fr{\alpha} & X \fd{p} \\
B \ar@{-->}[ru]^{g}\fr{\beta} & Y}
\]
there exists $g$ making both triangles commutative. \ed

\bd If $I$ is a set of morphisms of flows, the collection of 
morphisms of flows that satisfies the RLP with respect to any morphism
of $I$ is denoted by $I-inj$. Denote by $I-cof$ the collection of
morphisms of flows that satisfies the RLP with respect to any morphism
that satisfies the LLP with respect to any element of $I$. This is a
purely categorical fact that $I-cell\subset I-cof$. \ed

\bd A \textit{model structure} on a category $\C$ consists of
three subcategories of the category of morphisms $Map(\C)$ called weak
equivalences, cofibrations, and fibrations, and two functorial
factorizations $(\alpha,\beta)$ and $(\gamma,\delta)$ satisfying the
following properties:
\begin{enumerate}
\item (2-out-of-3) If $f$ and $g$ are morphisms of $\C$ such that $g\circ f$
is defined and two of $f$, $g$ and $g\circ f$ are weak
equivalences, then so is the third.
\item (Retracts) If $f$ and $g$ are morphisms of $\C$ such that $f$ is a retract
of $g$ and $g$ is a weak equivalence, cofibration, or fibration,
then so is $f$.
\item (Lifting) Define a map to be a trivial cofibration if it is both a cofibration
and a weak equivalence. Similarly, define a map to be a trivial
fibration if it is both a fibration and a weak equivalence. Then
trivial cofibrations have the LLP with respect to fibrations, and
cofibrations have the LLP with respect to trivial fibrations.
\item (Factorization) For any morphism $f$, $\alpha(f)$ is a cofibration,
$\beta(f)$ a trivial fibration, $\gamma(f)$ is a trivial
cofibration , and $\delta(f)$ is a fibration.
\end{enumerate}
\ed

\bd A \textit{model category} is a complete and cocomplete
category $\C$ together with a model structure on $\C$. \ed

\bth\label{sufficient} \cite{MR99h:55031} Let $\C$ be a complete
and cocomplete category. Let $\mathcal{W}$ be a subcategory of
$Map(\C)$. Let $I$ and $J$ be two sets of maps of $\C$. Then there
exists a structure of model category on $\C$ such that the fibrations
are exactly the arrows satisfying the RLP\index{right lifting
property} with respect to the arrows of $J$, such that the trivial
fibrations are exactly the arrows satisfying the RLP with respect to
the arrows of $I$, such that the weak equivalences are exactly the
arrows of $\mathcal{W}$ if the following conditions are satisfied:
\begin{enumerate}
\item The subcategory $\mathcal{W}$ has the 2-out-of-3 property and is
closed under retracts.
\item The domains of $I$ are small relative to $I$-cell.
\item The domains of $J$ are small relative to $J$-cell.
\item Any relative $J$-cell complex  is a weak equivalence and satisfies the
LLP with respect to any morphism satisfying the RLP\index{right
lifting property} with respect to the arrows of $I$. In other
terms, $J-cell\subset I-cof\cap \mathcal{W}$.
\item A morphism satisfies the RLP with respect to the morphisms of $I$ if and
only if it is a weak equivalence and it satisfies the RLP with respect
to the morphisms of $J$. In other terms, $I-inj=J-inj\cap
\mathcal{W}$.
\end{enumerate}
\eth

\bd If the conditions of Theorem~\ref{sufficient} are satisfied
for some model category $\C$, the set $I$ is the set of
\textit{generating cofibrations}, the set $J$ is the set of
\textit{generating trivial cofibrations} and one says that $\C$
is a \textit{cofibrantly generated model category}. \ed

The above conditions are satisfied for $\top$ if $\mathcal{W}$ is the
subcategory of weak homotopy equivalences, if $I$ is the set of
inclusion maps $\mathbf{S}^{n-1}\longrightarrow \mathbf{D}^n$ with
$\mathbf{S}^{-1}=\varnothing$ and for $n\geq 0$, and if $J$ is the set
of continuous maps $\mathbf{D}^n\longrightarrow [0,1]\p \mathbf{D}^n$
such that $x\mapsto (0,x)$ and for $n\geq 0$. The fibrations of the
model structure of $\top$ are usually called Serre fibration.

So far, we have proved:

\bth The category of flows $\dtop$ is complete and cocomplete. Moreover:
\begin{enumerate}
\item The subcategory $\mathcal{S}$ has the 2-out-of-3 property and is
closed under retracts.
\item The domains of $I^{gl}_+$ are small relative to $I^{gl}_+$-cell.
\item The domains of $J^{gl}$ are small relative to $J^{gl}$-cell.
\end{enumerate}
\eth

\section{Characterization of the fibrations of flows}\label{chafib}

\bd An element of $J^{gl}-inj$ is called a {\rm fibration}. A
fibration is {\rm trivial} if it is at the same time a weak S-homotopy
equivalence. \ed

\bp \label{ortho3} A morphism of flows
$f:X\longrightarrow Y$ satisfies the RLP with respect to
$\glob(U)\longrightarrow \glob(V)$ if and only if for any
$\alpha,\beta\in X^0$, $\P_{\alpha,\beta}X\longrightarrow
\P_{f(\alpha),f(\beta)}Y$ satisfies the RLP with respect to
$U\longrightarrow V$. \ep

\bpf Considering a commutative square of topological spaces
\[
\xymatrix{U \fr{}\fd{} & \P_{\alpha,\beta}X\fd{f} \\
V\fr{} \ar@{-->}[ru]^-{k_1}& \P_{f(\alpha),f(\beta)}Y}
\]
is equivalent to considering a commutative square of flows like
\[
\xymatrix{\glob(U) \fr{\begin{array}{c} 0\mapsto \alpha\\1\mapsto \beta\end{array}}\fd{} & X\fd{f} \\
\glob(V)\fr{} \ar@{-->}[ru]^-{k_2}& Y}
\]
The existence of $k_1$ making the first diagram commutative is
equivalent to the existence of $k_2$ making the second diagram
commutative. Hence the result.  \epf

\bp Let $f:X\longrightarrow Y$ be a morphism of flows. Then $f$ is a
fibration of flows if and only if $\P f:\P X\longrightarrow \P Y$ is
a Serre fibration of topological spaces. \ep

\bpf By Proposition~\ref{ortho3},
the morphism of flows $f$ is a fibration if and only if for any
$\alpha,\beta\in X^0$, the continuous map
$\P_{\alpha,\beta}X\longrightarrow
\P_{f(\alpha),f(\beta)}Y$ is a Serre fibration. But $\mathbf{D}^n$ and
$\mathbf{D}^n\p [0,1]$ are connected. Hence the result. \epf

\section{About the necessity of $R$ and $C$ as generating cofibrations}
\label{necessaire}

One cannot take as definition of a cofibration an element of
$I^{gl}-cof$. Indeed:

\bp There does not exist any cofibrantly generated model structure on $\dtop$
such that the generating set of cofibrations is $I^{gl}$, the
generating set of trivial cofibrations $J^{gl}$, and the class of weak
equivalences the one of weak S-homotopy equivalences. \ep

\bpf  If such a model structure existed, then all cofibrations
would be synchronized because any cofibration is a retract of an
element of $I^{gl}-cell$, because any element of $I^{gl}-cell$ is
synchronized, and at last because the retract of a synchronized
morphism of flows is synchronized. Since a trivial fibration is a weak
S-homotopy equivalence, then such morphism is in particular
synchronized. So all composites of the form $p\circ i$ where $p$ would
be a trivial fibration and $i$ a cofibration would be synchronized. So
a non-synchronized morphism of flows could never be equal to such
composite. \epf

\bp There does not exist any cofibrantly generated model structure on $\dtop$
such that the generating set of cofibrations is $I^{gl}\cup\{C\}$, the
generating set of trivial cofibrations $J^{gl}$, and the class of weak
equivalences the one of weak S-homotopy equivalences. \ep

\bpf Suppose that such a model structure exists. Consider a commutative
square
\[
\xymatrix{
A\fd{i} \fr{} & \{0,1\} \fd{R} \\
X \fr{} \ar@{-->}[ru]^{k}& \{0\}}
\]
where $i:A\longrightarrow X$ is an element of $I^{gl}\cup\{C\}$.
Since the path spaces of the flows $\{0,1\}$ and $\{0\}$ are empty,
then $\P A=\P X=\varnothing$. So $i=C$, $A=\varnothing$ and
$X=\{0\}$. Let $k(0)=0$. Then $k$ makes the diagram above
commutative. Therefore $R$ satisfies the RLP with respect to any
morphism of $I^{gl}\cup\{C\}$. So $R$ is a trivial fibration for this
model structure. Contradiction. \epf

\bp There does not exits any cofibrantly generated model structure on $\dtop$
such that the generating set of cofibrations is
$I^{gl}\cup\{R\}$, the generating set of trivial cofibrations
$J^{gl}$, and the class of weak equivalences the one of weak
S-homotopy equivalences. \ep

\bpf  If such a model structure existed, then all cofibrations
would restrict to an onto set map between the $0$-skeletons. So there
would not exist any cofibrant object since the initial flow is the
empty set. \epf

Hence the definition:

\bd An element of $I^{gl}_+-cof$ is called a
{\rm cofibration}.  A cofibration is {\rm trivial} if it is at the
same time a weak S-homotopy equivalence. \ed

\section{Pushout of $\glob(\de Z)\rightarrow \glob(Z)$ in $\dtop$}\label{expush}

Let $\de Z\longrightarrow Z$ be a continuous map.
Let us consider a diagram of flows as follows:
\[
\xymatrix{
\glob(\de Z) \fr{\phi}\fd{} & A  \fd{}\\
\glob(Z) \fr{} & \cocartesien X}
\]
The purpose of this short section is an explicit description of the
pushout $X$ in the category of flows.

Let us consider the set $\mathcal{M}$ of finite sequences
$\alpha_0\dots\alpha_p$ of elements of $A^0=X^0$ with $p\geq 1$
and such that, for any $i$, at least one of the two pairs
$(\alpha_i,\alpha_{i+1})$ and $(\alpha_{i+1},\alpha_{i+2})$ is
equal to $(\phi(0),\phi(1))$. Let us consider the pushout diagram
of topological spaces
\[
\xymatrix{
\de Z \fr{\phi}\fd{} & \P_{\phi(0),\phi(1)} A \fd{} \\
Z \fr{} & \cocartesien T }
\]

Let $Z_{\alpha,\beta}=\P_{\alpha,\beta}A$ if $(\alpha,\beta)\neq
(\phi(0),\phi(1))$ and let $Z_{\phi(0),\phi(1)}=T$. At last, for
any $\alpha_0\dots\alpha_p\in \mathcal{M}$, let
$[\alpha_0\dots\alpha_p]= Z_{\alpha_0,\alpha_1}\p
Z_{\alpha_1,\alpha_2}\p \dots \p Z_{\alpha_{p-1},\alpha_p}$.  And
$[\alpha_0\dots\alpha_p]_i$ denotes the same product as
$[\alpha_0\dots\alpha_p]$ except that
$(\alpha_i,\alpha_{i+1})=(\phi(0),\phi(1))$ and that the factor
$Z_{\alpha_{i},\alpha_{i+1}}=T$ is replaced by
$\P_{\phi(0),\phi(1)} A$.  We mean that in the product
$[\alpha_0\dots\alpha_p]_i$, the factor $\P_{\phi(0),\phi(1)} A$
appears exactly once. For instance, one has (with $\phi(0)\neq
\phi(1)$) \beas &&
[\alpha\phi(0)\phi(1)\phi(0)\phi(1)]=\P_{\alpha,\phi(0)}A\p T\p
\P_{\phi(1),\phi(0)}A\p T\\
&& [\alpha\phi(0)\phi(1)\phi(0)\phi(1)]_1=\P_{\alpha,\phi(0)}A\p
\P_{\phi(0),\phi(1)} A\p \P_{\phi(1),\phi(0)}A\p T\\
&& [\alpha\phi(0)\phi(1)\phi(0)\phi(1)]_3=\P_{\alpha,\phi(0)}A\p T\p
\P_{\phi(1),\phi(0)}A\p \P_{\phi(0),\phi(1)} A. \eeas The idea is that
in the products $[\alpha_0\dots\alpha_p]$, there are no possible
simplifications using the composition law of $A$. On the contrary,
exactly one simplification is possible using the composition law of
$A$ in the products $[\alpha_0\dots\alpha_p]_i$. For instance, with
the examples above, there exist continuous maps
\[[\alpha\phi(0)\phi(1)\phi(0)\phi(1)]_1\longrightarrow
[\alpha\phi(0)\phi(1)]\] and
\[[\alpha\phi(0)\phi(1)\phi(0)\phi(1)]_3\longrightarrow[\alpha\phi(0)\phi(1)\phi(1)]\]
induced by the composition law of $A$ and there exist
continuous maps
\[[\alpha\phi(0)\phi(1)\phi(0)\phi(1)]_1\longrightarrow [\alpha\phi(0)\phi(1)\phi(0)\phi(1)]\]
and
\[[\alpha\phi(0)\phi(1)\phi(0)\phi(1)]_3\longrightarrow [\alpha\phi(0)\phi(1)\phi(0)\phi(1)]\]
induced by the continuous map $\P_{\phi(0),\phi(1)}A\longrightarrow T$.

Let $\P_{\alpha,\beta}M$ be the colimit of the diagram of
topological spaces consisting of the topological spaces
$[\alpha_0\dots\alpha_p]$ and $[\alpha_0\dots\alpha_p]_i$ with
$\alpha_0=\alpha$ and $\alpha_p=\beta$ with the two kinds of maps
above defined. The composition law of $A$ and the free
concatenation obviously defines a continuous associative map
$\P_{\alpha,\beta}M\p \P_{\beta,\gamma}M\longrightarrow
\P_{\alpha,\gamma}M$.

\bp \label{pushexplicit}\label{pushexplicit0}
One has the pushout diagram of flows
\[
\xymatrix{
\glob(\de Z) \fr{\phi}\fd{} & A \fd{} \\
\glob(Z) \fr{} & \cocartesien M }
\]
\ep

\bpf Let us consider a commutative diagram like:
\[
\xymatrix{
\glob(\de Z)\fr{\phi} \fd{} & A \fd{}\ar@/^15pt/[rdd]^{\phi_1}&\\
\glob(Z) \ar@/_15pt/[rrd]_{\phi_2}\fr{} & M\ar@{-->}[rd]^{{h}} &\\
&&X}
\]
One has to prove that there exists $h$ making everything
commutative. We do not have any choice for the definition on the
$0$-skeleton: $h(\alpha)=\phi_1(\alpha)$. The diagram of flows above
gives a commutative diagram of topological spaces
\[
\xymatrix{
\de Z\fr{\phi} \fd{} & \P A \fd{}\ar@/^15pt/[rdd]^{\phi_1}&\\
Z \ar@/_15pt/[rrd]_{\phi_2}\fr{} & T\cocartesien\ar@{-->}[rd]^{{k}} &\\
&&\P X}
\]
By construction of $T$, there exists a continuous map
$k:T\longrightarrow \P_{h(\phi(0)),h(\phi(1))}X \subset \P X$ making
the diagram commutative.

Constructing a continuous map $\P M\longrightarrow \P X$ is equivalent
to constructing continuous maps
$[\alpha_0\dots\alpha_p]\longrightarrow
\P_{h(\alpha_0),h(\alpha_p)}X$ and
$[\alpha_0\dots\alpha_p]_i\longrightarrow
\P_{h(\alpha_0),h(\alpha_p)}X$ for any finite sequence
$\alpha_0\dots\alpha_p$ of $\mathcal{M}$ such that any diagram
like
\[
\xymatrix{[\alpha_0\dots\alpha_p]_i\fd{} \fr{}
&\P_{h(\alpha_0),h(\alpha_p)}X\\
[\alpha_0\dots\alpha_p]\ar@{->}[ru]&}
\]
or like
\[
\xymatrix{[\alpha_0\dots\alpha_p]_i\fd{} \fr{}
&\P_{h(\alpha_0),h(\alpha_p)}X\\
[\alpha_0\dots \widehat{\phi(0)\phi(1)}\dots\alpha_p]\ar@{->}[ru]&}
\]
is commutative. There are such obvious maps by considering the
continuous maps $Z_{\alpha,\beta}\longrightarrow
\P_{h(\alpha),h(\beta)}X$ and by composing with the composition
law of $X$. Hence the result. \epf

\bth\label{classiquen}\label{comp2} Suppose that one has the
pushout of flows
\[
\xymatrix{
\glob(P) \fr{\phi}\fd{} & A \fd{} \\
\glob(Q) \fr{} & \cocartesien X }
\]
where $P\longrightarrow Q$ is an inclusion of a deformation retract of
topological spaces. Then the continuous map $\P f:\P A\longrightarrow
\P X$ is a weak homotopy equivalence. \eth

\bpf Let us start with the diagram $\D=\D_0$ of topological spaces
constructed for Proposition~\ref{pushexplicit} calculating $\P X$.  We
are going to modify $\D$, by transfinite induction, in order to obtain
another diagram of topological spaces, whose colimit will still be
isomorphic to $\P X$ and such that all arrows will be inclusions of a
deformation retract.

We are going to add vertices and arrows to the diagram above in the
following way. For any configuration like
\[
\xymatrix{
[\alpha_0\dots\alpha_p]_i \fr{j}\fd{c} & [\alpha_0\dots\alpha_p]\\
[\alpha_0\dots\widehat{\alpha_i}\dots \alpha_p]}
\]
where $c$ is induced by the composition law of $A$ and $j$ is the
unique possible inclusion of a deformation retract, let us draw the
cocartesian square
\[
\xymatrix{
[\alpha_0\dots\alpha_p]_i \fr{j}\fd{c} & [\alpha_0\dots\alpha_p]\ar@{-->}[d]^-{d}\\
[\alpha_0\dots\widehat{\alpha_i}\dots \alpha_p]
\ar@{-->}[r]^-{k}& \cocartesien U}
\]
Notice that $k$ is an inclusion of a deformation retract because the
class of inclusions of a deformation retract is closed under pushout :
cf. \cite{MR99h:55031} for an elementary proof, or \cite{ruse} for a
model-categoric argument. Indeed, an inclusion of a deformation retract
is a trivial cofibration for the Str{\o}m model category of compactly
generated topological spaces. So the corresponding class is closed
under pushout because it coincides with the class of morphisms
satisfying the LLP with respect to any Hurewicz fibration.

One will say that the maps $j$ and $k$ are
\textit{orthogonal to the composition law of $A$} and that the
maps $c$ and $d$ are \textit{parallel to the composition law of
$A$}. Repeat the process for any configuration like
\[
\xymatrix{
U \fr{j}\fd{c} & V\\
W}
\]
where $j$ is orthogonal to the composition law of $A$ and $c$
parallel to the composition law of $A$ by completing the
configuration by a cocartesian square of topological spaces
\[
\xymatrix{
U \fr{j}\fd{c} & V\ar@{-->}[d]^-{d}\\
W \ar@{-->}[r]^-{k}& X \cocartesien}
\]
By induction, one will say that $k$ is orthogonal to the composition
law and that $d$ is parallel to the composition law.  Notice that, in
this diagram, any map which is orthogonal to the composition law is an
inclusion of a deformation retract of topological spaces. At each step
consisting of adding an object so that it creates a pushout square in
the diagram, one obtains a diagram $\D_{\lambda+1}$ from a diagram
$\D_\lambda$. There is a canonical continuous map $\liminj
\D_\lambda\longrightarrow \liminj \D_{\lambda+1}$ which is an homeomorphism.

Let us say that the topological spaces $[\alpha_0\dots\alpha_p]$ and
$[\alpha_0\dots\alpha_p]_i$ are of length $p$. By induction , one
defines the length of a topological space as being constant along the
arrows orthogonal to the composition law of $A$. The length is
strictly decreasing along the arrows parallel to the composition law
of $A$. Therefore the process stops after an, eventually, transfinite
number of steps. Moreover the only map which can starts from an
element of length $1$ is an arrow orthogonal to the composition law of
$A$. Therefore such a map is necessarily an inclusion of a deformation
retract of topological spaces.

Let us say that the process stops for $\lambda=\lambda_0$. For
any vertex $v$ of $\D_{\lambda_0}$, there exists an arrow
$v\longrightarrow w$ of $\D_{\lambda_0}$ with $w$ of length $1$.
Therefore the colimit of the diagram $\D_{\lambda_0}$ is
isomorphic to the colimit of the subdiagram of $\D_{\lambda_0}$
consisting of the vertex of length $1$.

The initial diagram $\D=\D_0$ has therefore the same colimit as a
diagram of topological spaces of the form a concatenation of
straight lines of the form
\[\P_{\alpha,\beta}A\longrightarrow M_1 \longrightarrow M_2 \longrightarrow \dots\]
where all arrows are inclusions of a deformation retract. Therefore
$\P A\longrightarrow \P X$ is a weak homotopy equivalence since any
inclusion of a deformation retract is a closed $T_1$ inclusion and a
weak homotopy equivalence and since any transfinite composition of such
maps is a weak homotopy equivalence (cf \cite{MR99h:55031} Lemma~2.4.5,
Corollary~2.4.6 and Lemma~2.4.8). \epf

\section{$J^{gl}-cell\subset I^{gl}_+-cof\cap \mathcal{S}$}\label{fin1}

\bp\label{z1} One has $J^{gl}-cell\subset I^{gl}-cof\cap
\mathcal{S}$. \ep

\bpf The continuous maps
$\mathbf{D}^n\iso \mathbf{D}^n\p\{0\}\longrightarrow \mathbf{D}^n\p
[0,1]$ are inclusions of a deformation retract for any $n\geq 0$. So
by Theorem~\ref{classiquen}, $J^{gl}-cell\subset\mathcal{S}$. The
class $I^{gl}-cof$ is closed under pushout and transfinite
composition. So it then suffices to prove that $J^{gl}\subset
I^{gl}-cof$. A morphism of flows $f:X\longrightarrow Y$ satisfies the
RLP with respect to $\glob(\mathbf{D}^n)\longrightarrow
\glob(\mathbf{D}^n\p [0,1])$ if and only if for any $\alpha,\beta\in
X^0$, $\P_{\alpha,\beta}X\longrightarrow
\P_{f(\alpha),f(\beta)}Y$ is a Serre fibration by Proposition~\ref{ortho3}.
But again by Proposition~\ref{ortho3}, for any element
$f:X\longrightarrow Y$ of $I^{gl}-inj$, for any $\alpha,\beta\in X^0$,
$\P_{\alpha,\beta}X\longrightarrow
\P_{f(\alpha),f(\beta)}Y$ is a trivial Serre fibration, 
so a Serre fibration. Hence the result.  \epf

\bp\label{ortho2} Let $f$ be a morphism of flows. Then the
following conditions are equivalent:
\begin{enumerate}
\item $f$ is synchronized
\item $f$ satisfies the RLP
with respect to $R:\{0,1\}\longrightarrow \{0\}$ and
$C:\varnothing\subset \{0\}$.
\end{enumerate}
 \ep

\bpf Let $f:X\longrightarrow Y$ that satisfies the RLP with respect
to $R:\{0,1\}\longrightarrow \{0\}$ and $C:\{0\}\subset \{0,1\}$. Let
us suppose that $f(a)=f(b)$ for some $a,b\in X^0$. Then consider
the commutative diagram
\[
\xymatrix{
\{0,1\} \ar@{->}[rr]^-{0\mapsto a,1\mapsto b} \fd{R} && X \fd{f}\\
\{0\} \ar@{-->}[rru]^{g}\ar@{->}[rr]^-{0\mapsto f(a)} && Y}
\]
By hypothesis, there exists $g$ making both triangles
commutative. So $b=g\circ R(1)=g(0)=g\circ R(0)=a$. So $f$
induces a one-to-one map on the $0$-skeletons. Now take $a\in
Y^0$. Then consider the commutative diagram
\[
\xymatrix{
\varnothing \ar@{->}[rr]^-{} \fd{C} && X \fd{f}\\
\{0\} \ar@{-->}[rru]^{g}\ar@{->}[rr]^-{0\mapsto a} && Y}
\]
By hypothesis, there exists $g$ making both triangles
commutative. Then $a=f(g(1))$. So $f$ induces an onto map on the
$0$-skeletons. Therefore condition~2 implies condition~1.  Conversely,
if $f$ is synchronized, let $(f^0)^{-1}:Y^0\longrightarrow X^0$ be the
inverse of the restriction $f^0$ of $f$ to the $0$-skeleton. Consider
a commutative diagram like
\[
\xymatrix{
A \ar@{->}[rr]^-{\alpha} \fd{u} && X \fd{f}\\
B \ar@{-->}[rru]^{g}\ar@{->}[rr]^-{\beta} && Y}
\]
where $A$ and $B$ are two flows such that $A=A^0$, $B=B^0$ and where
$u$ is any set map. Then $g=(f^0)^{-1}\circ \beta$ makes both
triangles commutative. Indeed $f\circ (f^0)^{-1}\circ
\beta=\beta$ and $(f^0)^{-1}\circ \beta\circ u=(f^0)^{-1}\circ
f\circ \alpha=\alpha$. \epf

\begin{cor} $J^{gl}-cell\subset I^{gl}_+-cof\cap \mathcal{S}$.
\end{cor}

\bpf The elements of $I^{gl}-cof$ satisfies the LLP with respect
to any element of $I^{gl}-inj$. So in particular, the elements of
$I^{gl}-cof$ satisfies the LLP with respect to any synchronized
element of $I^{gl}-inj$. But a synchronized element of $I^{gl}-inj$ is
precisely an element of $I^{gl}_+-inj$ by
Proposition~\ref{ortho2}. Therefore $I^{gl}-cof\subset
I^{gl}_+-cof$. \epf

\section{$I^{gl}_+-inj= J^{gl}-inj\cap
\mathcal{S}$}\label{fin2}

\bp \label{z2} Any morphism of $I^{gl}_+-inj$ is a trivial
fibration. In other terms, $I^{gl}_+-inj\subset J^{gl}-inj\cap
\mathcal{S}$. \ep

\bpf Let $f:X\longrightarrow Y$ be a morphism of flows with $f\in I^{gl}_+-inj$.
Then $f$ is synchronized by Proposition~\ref{ortho2}. By
Proposition~\ref{ortho3}, for any $\alpha,\beta\in X^0$, $\P_{\alpha,\beta}X\longrightarrow
\P_{f(\alpha),f(\beta)}Y$ is a trivial Serre fibration. So $f$ is a weak S-homotopy
equivalence. And again by Proposition~\ref{ortho3}, this implies that
$f$ satisfies the RLP with respect to $J^{gl}$. Hence the result. \epf

\bp\label{z3} Any trivial fibration is in $I^{gl}_+-inj$. In other
terms, $J^{gl}-inj\cap \mathcal{S}\subset I^{gl}_+-inj$.\ep

\bpf Let $f$ be a trivial fibration. By Proposition~\ref{ortho3},
for any $\alpha,\beta\in X^0$, the continuous map
$\P_{\alpha,\beta}X\longrightarrow \P_{f(\alpha),f(\beta)}Y$ is a
fibration. But $f\in \mathcal{S}$. Therefore the fibrations
$\P_{\alpha,\beta}X\longrightarrow \P_{f(\alpha),f(\beta)}Y$ are
trivial. So by Proposition~\ref{ortho3}, $f$ satisfies the RLP with
respect to $I^{gl}$. Since $f$ is also synchronized, then $f$
satisfies the RLP with respect to $R$ and $C$ as well.  \epf

\section{The model structure of $\dtop$}\label{modelflow}

\begin{cor}
\label{c} The category of flows together with the weak S-homotopy equivalences,
the cofibrations and the fibrations is a model\index{model category}
category. The cofibrations are the retracts of the elements of
$I^{gl}_+-cell$. Moreover, any flow is fibrant.
\end{cor}

\bpf The first part of the statement is a consequence of
Proposition~\ref{compact}, Proposition~\ref{z1}, Proposition~\ref{z2},
Proposition~\ref{z3} and Theorem~\ref{sufficient}. It remains to prove
that any flow is fibrant.  Let $X$ be a flow. Let $\mathbf{1}$ be the
flow such that $\mathbf{1}^0=\{0\}$ and $\P\mathbf{1}=\{1\}$. Then
$\mathbf{1}$ is a terminal object of $\dtop$. Consider a commutative
diagram like
\[
\xymatrix{
\glob(\mathbf{D}^n) \fd{i_0} \fr{\alpha} & X \fd{}\\
\glob([0,1]\p \mathbf{D}^n)\ar@{-->}[ru]^{g} \fr{\beta}& \mathbf{1}}
\]
Let $g(0)=\alpha(0)$, $g(1)=\alpha(1)$ and $g(t,z)=\alpha(z)$ for
any $(t,z)\in [0,1]\p \mathbf{D}^n$. Then $g$ makes both triangles
commutative. \epf

\begin{cor} Any cofibration for this model structure induces a closed
inclusion between path spaces. \end{cor}

\bpf This is a consequence of Corollary~\ref{pushout-Scof} and of
Corollary~\ref{c}. \epf

\section{S-homotopy and the model structure of $\dtop$}\label{comparison}

In any model category, the canonical morphism $X\sqcup
X\longrightarrow X$ factors as a cofibration $X\sqcup X\longrightarrow
I(X)$ and a trivial fibration $I(X)\longrightarrow X$. One then says
that two morphisms $f$ and $g$ from $X$ to $Y$ are left homotopy
equivalent (this situation being denoted by $f\sim_l g$) if and only
if there exists a morphism $I(X)\longrightarrow Y$ such that the
composite $X\sqcup X\longrightarrow I(X)\longrightarrow Y$ is exactly
$f\sqcup g$. On cofibrant and fibrant objects, the left homotopy is an
equivalence relation simply called homotopy. Then one can say that two
cofibrant and fibrant flows $X$ and $Y$ are left homotopy equivalent
(this situation being denoted by $X\sim_l Y$) if and only if there
exists a morphism of flows $f:X\longrightarrow Y$ and a morphism of
flows $g:Y\longrightarrow X$ such that $f\circ g\sim_l \Id_Y$ and
$g\circ f\sim_l \Id_X$.

\bth\label{pareil} Two cofibrant flows are left homotopy
equivalent if and only if they are S-homotopy equivalent. \eth

The similar fact is trivial in $\top$ because for any cofibrant
topological space $X$, the continuous map $X\sqcup X\longrightarrow
[0,1]\p X$ sending one copy of $X$ to $\{0\}\p X$ and the other one to
$\{1\}\p X$ is a relative $I$-cell complex, and therefore a
cofibration for the model structure of $\top$, and the continuous
projection map $[0,1]\p X\longrightarrow X$ is a fibration. A similar
situation does not hold in the framework of flows.

\bp There exists a cofibrant flow $X$ such that the canonical
morphism of flows $[0,1]\boxtimes X\longrightarrow X$ such that
$t\boxtimes x\mapsto x$ is not a fibration for the model
structure of $\dtop$. \ep

\bpf Let $X^0$ be the three-element set
$\{\alpha,\beta,\gamma\}$. Let $\P_{\alpha,\beta}X=\{u\}$,
$\P_{\beta,\gamma}X=\{v\}$, and $\P_{\alpha,\gamma}X=\mathbf{D}^1$ with the
relation $1=u*v$. Consider the commutative diagram
\[
\xymatrix{
\glob(\mathbf{S}^0)\fd{}\fr{\alpha}& [0,1]\boxtimes X\fd{}\\
\glob(\mathbf{D}^1)\fr{\beta} & X}
\]
with $\alpha(-1)=0\boxtimes -1$, $\alpha(1)=(0\boxtimes u)*(1\boxtimes
v)$ and $\beta(z)=z$ for $z\in \mathbf{D}^1$. Suppose that there
exists $g:\glob(\mathbf{D}^1)\longrightarrow [0,1]\boxtimes X$ making the
above diagram commutative. For $z\in \mathbf{D}^1\backslash\{1\}$,
then the execution path $t\boxtimes z$ of $[0,1]\boxtimes X$ is
composable with nothing by construction. So for such $z$,
$g(z)=\phi(z)\boxtimes z$ for some continuous map
$\phi:\mathbf{D}^1\backslash\{1\}\longrightarrow [0,1]$. Then
$n\mapsto\phi(1-1/(n+1))$ is a sequence of $[0,1]$ and so contains a
subsequence converging to some $t_0\in [0,1]$.  Then $(0\boxtimes
u)*(1\boxtimes v)=t_0\boxtimes 1$, which contredicts the explicit
description of $[0,1]\boxtimes X$ of
Corollary~\ref{disparition1}. \epf

\bd Let $X$ be a flow. Then the flow $\square X$ is defined by the
cocartesian diagram
\[
\xymatrix{
X^0\fd{i} \fr{i} & X \fd{} \\
X\fr{} & \cocartesien \square X}
\]
where $i:X^0\longrightarrow X$ is the canonical inclusion. This flow
is called the \textit{square of $X$}. \ed

\bp\label{z4} For any flow $X$, the canonical morphism of flows
$k_X:X\sqcup X\longrightarrow \square X$ is a cofibration. \ep

\bpf This map is indeed an (eventually transfinite) composition of
pushouts of $R:\{0,1\}\longrightarrow \{0\}$, so an element of
$I^{gl}_+-cell\subset I^{gl}_+-cof$. \epf

\bp Let $X$ be a cofibrant flow. Then the canonical morphism of
flows $j_X:X\sqcup X\longrightarrow [0,1]\boxtimes X$ induced by the
inclusions $X\iso \{0\}\boxtimes X\subset [0,1]\boxtimes X$ and
$X\iso \{1\}\boxtimes X\subset [0,1]\boxtimes X$ is a
cofibration. \ep

\bpf The morphism $j_X$ factors as $j_X=\ell_X\circ k_X$. Using
Proposition~\ref{z4}, it suffices to prove that $\ell_X:\square
X\longrightarrow [0,1]\boxtimes X$ is a cofibration. Both functors
$X\mapsto \square X$ and $X\mapsto [0,1]\boxtimes X$ commute with
colimits and a colimit of cofibrations is a cofibration. So it
suffices to prove that for any CW-complex $Z$,
$\ell_{\glob(Z)}:\square {\glob(Z)}\longrightarrow [0,1]\boxtimes
{\glob(Z)}$ is a cofibration. But $\square {\glob(Z)}\iso
\glob(Z\sqcup Z)$ and $[0,1]\boxtimes {\glob(Z)}\iso
\glob([0,1]\p Z)$. Since $Z\sqcup Z\longrightarrow [0,1]\p Z$ is a
cofibration in $\top$, then it is a retract of an element of
$I-cell$. So $\ell_{\glob(Z)}$ is a retract of an element of
$I^{gl}-cell$. Therefore $\ell_{\glob(Z)}$ is a cofibration of
flows. \epf

\bp\label{part} Let $f:X\longrightarrow Y$ be a morphism of flows. If
$f$ is a S-homotopy\index{S-homotopy} equivalence, then it is
synchronized and for any $\alpha,\beta\in X^0$, the continuous map
$\P_{\alpha,\beta}f:\P_{\alpha,\beta}X\longrightarrow
\P_{f\left(\alpha\right),f\left(\beta\right)}Y$ is an homotopy
equivalence. \ep

\bpf Let $g:Y\longrightarrow X$ be a S-homotopic inverse of $f$. Let
$F:[0,1]\boxtimes X\longrightarrow X$ be a
S-homotopy\index{S-homotopy} from $g\circ f$ to $\Id_X$. Let
$G:[0,1]\boxtimes Y\longrightarrow Y$ be a
S-homotopy\index{S-homotopy} from $f\circ g$ to $\Id_Y$. Let
$\alpha,\beta\in X^0$. Then the composite
\[F':\glob\left([0,1]\p \P_{\alpha,\beta}X\right)\iso
[0,1]\boxtimes \glob\left(\P_{\alpha,\beta}X\right) \longrightarrow
[0,1]\boxtimes X\longrightarrow X\] defined by
\begin{itemize}
\item[] $F'\left(0\right)=f\left(\alpha\right)$ and $F'(1)=f\left(\beta\right)$
\item[] $F'\left(t,x\right)=F\left(t\boxtimes x\right)$ for $\left(t,x\right)\in [0,1]\p
\P_{\alpha,\beta}X$
\end{itemize}
yields an homotopy from $\P_{\alpha,\beta}g\circ\P_{\alpha,\beta}f$ to
$\Id_{\P_{\alpha,\beta}X}$. In the same way, one constructs from $G$
an homotopy from $\P_{\alpha,\beta}f\circ\P_{\alpha,\beta}g$ to
$\Id_{\P_{\alpha,\beta}Y}$. \epf

\bpf[Proof of Theorem~\ref{pareil}] The canonical morphism of
flows $[0,1]\boxtimes X\longrightarrow X$ factors as a cofibration
$i_X:[0,1]\boxtimes X\longrightarrow I(X)$ followed by a trivial
fibration $p_X:I(X)\longrightarrow X$. So the canonical morphism
$X\sqcup X\longrightarrow X$ factors as a cofibration $i_X\circ
j_X:X\sqcup X\longrightarrow I(X)$ followed by a trivial fibration
$p_X:I(X)\longrightarrow X$.

Now let $X$ and $Y$ be two S-homotopy equivalent cofibrant flows.
Let $f:X\longrightarrow Y$ be a S-homotopy equivalence between $X$
and $Y$. Then $f$ is a weak S-homotopy equivalence by
Proposition~\ref{part}. So $X$ and $Y$ are left-homotopy
equivalent. Reciprocally, let $X$ and $Y$ be two left-homotopy
equivalent cofibrant flows. Therefore there exists a morphism of
flows $f:X\longrightarrow Y$ and a morphism of flows $g:Y\longrightarrow
X$ such that $f\circ g\sim_l \Id_Y$ and $g\circ f\sim_l \Id_X$.
Then there exists morphisms of flows $H_X:I(X)\longrightarrow X$ and
$H_Y:I(Y)\longrightarrow Y$ such that $H_X\circ i_X\circ j_X=(g\circ
f)\sqcup \Id_X$ and $H_Y\circ i_Y\circ j_Y=(f\circ g)\sqcup
\Id_Y$. So $H_X\circ i_X$ is a S-homotopy between $g\circ f$ and
$\Id_X$ and $H_Y\circ i_Y$ is a S-homotopy between $f\circ g$ and
$\Id_Y$. Therefore $X$ and $Y$ are S-homotopy equivalent. \epf

One obtains finally the following theorem:

\bth \label{wW} There exists a structure of model\index{model
category} category on $\dtop$ such that
\begin{enumerate}
\item The weak equivalences are the weak S-homotopy equivalences.
\item The fibrations are the morphisms of flows that satisfy the RLP
with respect to the morphisms of flows
$\glob(\mathbf{D}^{n})\longrightarrow
\glob([0,1]\p \mathbf{D}^n)$ induced by the maps $x\mapsto (0,x)$.
\item The cofibrations are the morphisms of flows that satisfy the LLP
with respect to any morphism of flows that satisfies the RLP with
respect to the morphisms of flows
$\glob(\mathbf{S}^{n-1})\longrightarrow \glob(\mathbf{D}^n)$ induced
by the inclusions $\mathbf{S}^{n-1}\subset \mathbf{D}^n$ and with
respect to $R:\{0,1\}\longrightarrow \{0\}$ and $C:\varnothing\subset
\{0\}$.
\item Any flow is fibrant.
\item The fibration are the morphism of flows inducing a Serre fibration
of topological spaces between path spaces.
\item Two cofibrant flows are homotopy equivalent for this model
structure if and only if they are S-homotopy equivalent.
\end{enumerate}
\eth

\begin{cor} Let $X$ and $Y$ be two cofibrant flows. Let
$f:X\longrightarrow Y$ be a synchronized morphism of flows. Then the
following conditions are equivalent:
\begin{enumerate}
\item for any $\alpha,\beta\in X^0$, the continuous map
$\P_{\alpha,\beta}X\longrightarrow
\P_{f\left(\alpha\right),f\left(\beta\right)}Y$ is a weak
homotopy equivalence
\item for any $\alpha,\beta\in X^0$, the continuous map
$\P_{\alpha,\beta}X\longrightarrow
\P_{f\left(\alpha\right),f\left(\beta\right)}Y$ is homotopy
equivalence
\item $f$ is a weak S-homotopy equivalence
\item $f$ is a S-homotopy equivalence.
\end{enumerate}
\end{cor}

\begin{qu} How to find two flows $X$ and $Y$ (necessarily
not cofibrant) and a synchronized morphism of flows
$f:X\longrightarrow Y$ which is not a S-homotopy\index{S-homotopy}
equivalence and such that for any $\alpha,\beta\in X^0$, $f$
induces an homotopy equivalence from $\P_{\alpha,\beta}X$ to
$\P_{f(\alpha),f(\beta)}Y$.
\end{qu}

\section{Why no identity maps in the notion of flow ?}

There exist several reasons. Here is one of them. The section ``Why
non-contract\-ing maps ?'' of \cite{diCW} is also related to this
question. A similar phenomenon appears in the construction of the
``corner homology'' of an $\omega$-category in \cite{Gau}
(cf. Proposition~4.2 of the latter paper).

Let $X$ be a flow. Let us consider the topological space $\P^-X$ which
is solution of the following universal problem: there exists a
continuous map $h^-:\P X\longrightarrow \P^-X$ such that $h(x*y)=h(x)$
and any continuous map $f:\P X\longrightarrow Y$ such that
$f(x*y)=f(x)$ factors uniquely as a composite $\overline{f}\circ h^-$
for a unique continuous map $\overline{f}:\P^-X\longrightarrow Y$. And
let us consider the topological space $\P^+X$ which is solution of the
following universal problem: there exists a continuous map $h^+:\P
X\longrightarrow \P^+X$ such that $h(x*y)=h(y)$ and any continuous map
$f:\P X\longrightarrow Y$ such that $f(x*y)=f(y)$ factors uniquely as
a composite $\overline{f}\circ h^+$ for a unique continuous map
$\overline{f}:\P^+X\longrightarrow Y$. The space $\P^-X$ is called the
\textit{branching space} of $X$ and the space $\P^+X$ is called the
\textit{merging space} of $X$. Both mappings
$\P^-:\dtop\longrightarrow \top$ and $\P^+:\dtop\longrightarrow \top$
are crucial for the definition of T-homotopy (cf. \cite{model2}
\cite{hobr}).

Suppose now that a flow $X$ is a small category enriched over the
category of compactly generated topological spaces $\top$, that is we
suppose that there exists an additional continuous map
$i:X^0\longrightarrow \P X$ with $s(i(\alpha))=\alpha$ and
$t(i(\alpha))=\alpha$ for any $\alpha\in X^0$. Then for any $x\in \P
X$, we would have $x=s(x)*x$ and $x=x*t(x)$. So both topological
spaces $\P^-X$ and $\P^+X$ would be discrete. Therefore, in such a
setting, the correct definition would be for $\P^-X$ (resp. $\P^+X$)
the quotient of $\P X\backslash i(X^0)$ by the identifications $x=x*y$
(resp. $y=x*y$). But with such a definition, the mappings $X\mapsto
\P^-X$ and $X\mapsto \P^+X$ cannot be functorial anymore.

\section{Concluding discussion}

If $Z$ is a cofibrant topological space, then $\glob(Z)$ is a
cofibrant flow. Let us denote by $\top_c$ the full and faithful
subcategory of cofibrant topological spaces. Let us denote by
$\dtop_c$ the full and faithful subcategory of cofibrant flows.  Then
one has the commutative diagram of functors
\[\xymatrix{\top_c \fr{}\fd{\glob(-)}& \top\fd{\glob(-)}\\
\dtop_c\fr{} & \dtop}\] which becomes the commutative diagram of
functors
\[\xymatrix{\top_c[\mathcal{SH}^{-1}] \fr{\simeq}\ar@{^{(}->}[d]^-{\glob(-)}& \top[\mathcal{W}^{-1}]\ar@{^{(}->}[d]^-{\glob(-)}
\\
\dtop_c[\mathcal{SH}^{-1}]\fr{\simeq} &
\dtop[\mathcal{S}^{-1}]}\] where $\mathcal{SH}$ is the class
of homotopy equivalences (of topological spaces or of flows),
$\mathcal{W}$ the class of weak homotopy equivalences of topological
spaces, and at last $\mathcal{S}$ the class of weak S-homotopy
equivalences of flows. Both horizontal arrows of the latter diagram
are equivalence of categories. The notation $\C[\mathcal{X}^{-1}]$
means of course the localization of the category $\C$ with respect to
the class of morphisms $\mathcal{X}$.

\end{document}